\documentclass[12pt]{amsart}
\usepackage{amsfonts}
\usepackage{amsmath}
\usepackage{amssymb}
\usepackage{url,epsfig,eqngroup}
\usepackage{color}
\usepackage{graphicx}
\usepackage{subfigure}
\usepackage{cancel}

\newcommand{\cred}[1]{{\color{red} #1}}

\def\real{\mathbb{R}}
\def\cmplx{\mathbb{C}}
\def\ganz{\mathbb{Z}}

\def\eps{\varepsilon}

\newcommand\cF{\mathcal{F}}

\newcommand\cO{\mathcal{O}}

\renewcommand{\Re}{\operatorname{Re}}
\renewcommand{\Im}{\operatorname{Im}}

\newcommand{\smalf}{\par\smallskip\noindent}
\newcommand{\medlf}{\par\medskip\noindent}
\newcommand{\biglf}{\par\bigskip\noindent}

\newcommand{\be}{\begin{eqnarray}}
\newcommand{\ben}{\begin{eqnarray*}}
\newcommand{\en}{\end{eqnarray}}
\newcommand{\enn}{\end{eqnarray*}}
\newcommand{\Z}{{\mathbb Z}}
\newcommand{\N}{{\mathbb N}}

\newcommand{\R}{{\mathbb R}}
\newcommand{\G}{{\Gamma}}

\newtheorem{theorem}{Theorem}[section]
\newtheorem{definition}[theorem]{Definition} 
\newtheorem{lemma}[theorem]{Lemma}
\newtheorem{corollary}[theorem]{Corollary}

\newtheorem{remark}[theorem]{Remark}

\newtheorem{assumption}[theorem]{Assumption}

\definecolor{rot}{rgb}{0,0,0}

\definecolor{red}{rgb}{0,0,0}
\newcommand{\tr}{\textcolor{red}}

\definecolor{red1}{rgb}{0,0,0}
\newcommand{\trr}{\textcolor{red1}}

\definecolor{red2}{rgb}{0,0,0}
\newcommand{\tcrr}{\textcolor{red2}}

\begin{document}

\title{Time-harmonic scattering by \tcrr{locally perturbed} periodic structures with Dirichlet and Neumann boundary 
conditions}

\author{Guanghui Hu}
\address{Guanghui Hu: School of Mathematical Sciences and LPMC\\
Nankai University \\
Tianjin 300071, China}
\email{ghhu@nankai.edu.cn}
\author{Andreas Kirsch}
\address{Andreas Kirsch: Department of Mathematics \\
Karlsruhe Institute of Technology (KIT) \\
76131 Karlsruhe, Germany}
\email{andreas.kirsch@kit.edu}

\date{\today}

\begin{abstract}
The paper is concerned with well-posedness of TE and TM polarizations of time-harmonic 
electromagnetic scattering by perfectly conducting periodic surfaces and periodically arrayed 
obstacles \tcrr{with local perturbations}. \tcrr{The classical Rayleigh Expansion radiation condition does not always lead to well-posedness of the Helmholtz equation even in unperturbed periodic structures.}
We 
propose \trr{two equivalent radiation conditions} to characterize the radiating behavior of time-harmonic wave fields incited by a source term in an open waveguide under impenetrable boundary conditions. With \trr{these open waveguide radiation conditions}, uniqueness and existence of time-harmonic scattering by incoming point 
source waves, plane waves and \tr{surface waves} from locally perturbed periodic structures are 
established under either the Dirichlet or Neumann boundary condition.
 \tr{A Dirichlet-to-Neumann operator without using the Green's function is constructed for proving well-posedness of perturbed scattering problems.}

\vspace{.2in} 

\end{abstract}

\maketitle

\section{Introduction}
\begin{figure}[htb]
  \centering
  \subfigure[A Lipschitz periodic curve]{
  \includegraphics[width=0.5\textwidth]{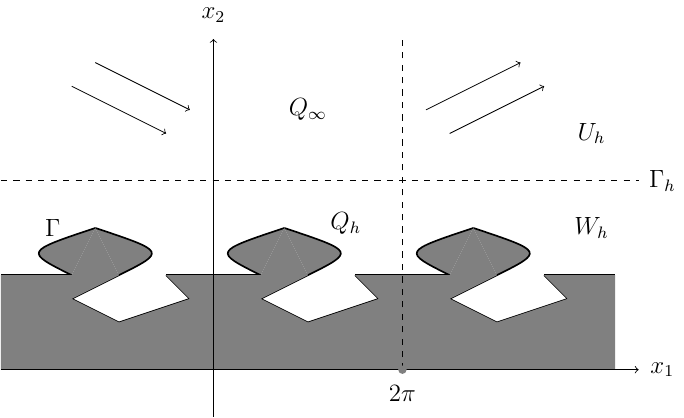}}\\
  \subfigure[Periodically arrayed obstacles]{
  \includegraphics[width=0.5\textwidth]{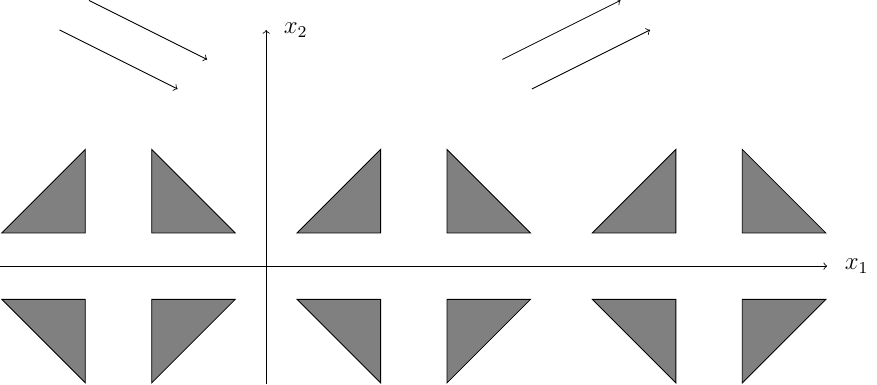}}
  \caption{Illustration of wave scattering from (a)  a perfectly reflecting periodic curve and 
  (b) perfectly conducting obstacles. 
\tcrr{Guided waves might exist in (a)-(b), leading to difficulties in establishing well-posedness of the scattering problem with the classical Rayleigh Expansion radiation condition \eqref{exc:b}. }}\label{f12}
\end{figure}

The electromagnetic scattering theory in periodic structures has many applications in 
micro-optics, radar imaging and non-destructive testing. We refer to \cite{P1980} for 
historical remarks and details of these applications. As a standard model, we consider a 
time-harmonic electromagnetic plane wave incident onto a  perfectly reflecting periodic 
surface or periodically arrayed conducting obstacles which remain invariant in the 
$x_3$-direction. Without loss of generality the direction of periodicity is supposed to 
be  $x_1$ and the arrayed obstacles lie in a layer of finite height in the $x_2$-direction.
We consider both the TE polarization case where the electric field is transversal to the 
$ox_1x_2$-plane by assuming $E(x)=(0,0,u(x_1,x_2))$ and the TM polarization case where 
the magnetic field is transversal to the $ox_1x_2$-plane by assuming $H(x)=(0,0,u(x_1,x_2))$. 
The background medium above the periodic surface or in the exterior of the periodically 
arrayed obstacles is supposed to be homogeneous and isotropic. The time-harmonic Maxwell's 
equations for $(E(x), H(x))$ will be reduced to the scalar Helmholtz equation for $u(x_1, x_2)$ 
over the $ox_1x_2$-plane together with the Dirichlet/Neumann boundary condition in TE/TM case 
and with proper radiation conditions as $|x_2|\rightarrow\infty$; see Figure \ref{f12} (a) 
and (b) for illustration of the scattering problems. 
\smalf
In periodic structures, a frequently used radiation condition is the so-called quasi-periodic 
Rayleigh expansion (see \eqref{exc:b}), which was firstly used by Lord Rayleigh in 1907 
\cite{R07}  for plane wave incidence. The Rayleigh expansion consists of a finite number of 
plane waves and infinitely many evanescent waves. However, such a radiation condition does 
not always lead to uniqueness of solutions for all frequencies due to the presence of 
evanescent/surface waves propagating along the unbounded periodic curve, or due to the 
existence of guided waves propagating between the arrayed obstacles, both of them decaying 
exponentially in $x_2$. Examples of surface waves for unbounded periodic curves of Dirichlet 
kind were constructed in \cite{G00} where the reflecting curve is not a graph and in \cite{KN02}
under the Neumann boundary condition. We also refer to \cite{BBS94} for non-uniqueness examples 
of solutions incited by periodically arrayed obstacles immersed in a dielectric layer. On the 
other hand, it is well known that surface waves do not exist if a Dirichlet periodic curve is 
given by the graph of some function \tcrr{or satisfies the geometrical condition \eqref{GC}}; see \cite{SM05, EY02, K93} for different regularity and 
geometry assumptions made on the reflecting curve. We also mention that the Rayleigh expansion 
condition does not apply to incoming source waves given by the fundamental solution of the 
Helmholtz equation and does not hold for scattering by compactly supported source terms. In 
these cases the incident waves lose the quasi-periodicity in $x_1$. It was firstly
discussed in \cite{CRZ98} that the radiated field should satisfy a Sommerfeld-type radiation 
condition and was recently proved \tcrr{in \cite{HWR} for Dirichlet rough surfaces given by graphs and in \cite{K22} for periodic inhomogeneous layers}. Hence, the radiating behavior of wave 
fields in periodic structures also depends on the type of incident waves.
\tcrr{To sum up, precise and sharp radiation conditions are still needed in order to mathematically interpret the radiating behavior of time-harmonic wave fields in  periodic structures, in particular for non-quasiperiodic incoming waves or when guided waves exist.} 
\smalf
In recent years, a new radiation condition has been derived from the limiting absorption 
principle for scattering by  layered periodic media in $\mathbb{R}^2$ and by periodic tubes 
in $\mathbb{R}^3$; see \cite{TF, K19-1, K19-2, KL18, K22, KL-MMAS}. Such a radiation condition
turns out to be equivalent to the radiation condition based on dispersion curves for closed 
periodic wave guides (see e.g., \cite{FJ16} and \cite[Remark 2.4]{K22}). By this new radiation 
condition, the diffracted fields \tcrr{caused by a compactly supported source term or a local defect} can be decomposed into the sum of a radiating part and a 
propagating (guided) part. The former decays as $|x_1|^{-3/2}$ in the horizontal direction 
$x_1$ and decays as $|x|^{-1/2}$ in the radical direction, whereas the latter is a finite 
number of quasi-periodic left-going and right-going evanescent modes which decay exponentially 
in the vertical direction $x_2$ (\cite{K22}). 
Moreover, this new radiation 
condition is stronger than the angular spectrum representation \cite{SM05} and the upward 
propagating radiation condition \cite{CZ98} for rough surface scattering problems. It can also 
be used for proving well-posedness of scattering by locally perturbed inhomogeneous layers in 
the presence of guided waves; see \cite{TF, KL18, K22}. \smalf
The aim of this paper is to investigate well-posedness of time-harmonic scattering by locally 
perturbed periodic curves and periodically arrayed obstacles of Dirichlet and Neumann kinds.  
\tcrr{The main results of this paper are summarized as follows.
\begin{itemize}
\item[(i)] Propose two equivalent radiation conditions to prove uniqueness of weak solutions for periodic Lipschitz interfaces with local perturbations. The first radiation condition was adapted from \cite{KL18, K22} for characterizing left-going and right-going evanescent waves of the propagating part of wave fields.  It is referred to as the open waveguide radiation condition, in comparision with the closed waveguide radiation condition of \cite{FJ16}.
The second radiation condition, which modifies the asymptotic behavior of radiating part of the first one,  was motivated by the Sommerfeld radiation condition justified in \cite{HWR} and \cite[Section 6]{K22} for point source waves. The second radiation condition extends the well-posedness result of \cite{HWR} to general periodic Lipschitz curves of Dirichlet or Neumann kind, in particular when guides waves are present. Since the decaying condition of Sommerfeld type contains more information on the radiating part, the second radiation condition yields a simplified proof of the uniqueness; see Theorem \ref{Th2.15}. 
\item[(ii)] Existence of solutions for incoming plane waves, surfaces waves and point source waves in a locally perturbed periodic structure under a priori assumptions (Sections \ref{s-pert}). 
Unlike the scattering by inhomogeneous periodic layers with local perturbations 
\cite{KL18, K22, TF}, there is no analogue of the Lippmann-Schwinger integral equation under 
the Dirichlet and Neumann boundary conditions. This leads to difficulties in the analysis of 
wave scattering from perfectly reflecting periodic curves with local perturbations. Our idea 
is to reduce the scattering problem to a bounded domain enclosing the perturbed part by 
constructing the DtN operator.
For this purpose, we construct
a Dirichlet-to-Neumann operator without using the Green's function for proving well-posedness of the perturbed scattering problem.
\end{itemize}
}

\biglf

The remaining part of the paper is organized as follows. We 
first consider the perturbed/unperturbed scattering problem due to a compact source term. In
Section \ref{s-intro}, we describe an open waveguide radiation condition \trr{and its equivalent version, and use them} to prove the uniqueness results. In 
comparison with the results for layered media \cite{KL18, K22}, a more general transmission 
problem and the scattering by exponentially decaying source terms without a compact support
will be investigated in the unperturbed periodic domain (see Theorems \ref{t-2} and \ref{t-3}).
In Section \ref{s-pert}, we prove well-posedness results for incoming point source waves, plane 
waves as well as surface waves in the perturbed setting. Finally, concluding remarks will made 
in Section \ref{sec:N} on how to carry out the analysis for unbounded periodic Dirichlet curves 
to Neumann curves and to periodically arrayed obstacles with boundary conditions.

\section{Scattering by Dirichlet periodic curves with local perturbations: radiation 
condition and uniqueness}
\label{s-intro}
\subsection{Notations}
Let $D\subset\R^2$ be a $2\pi$-periodic domain with respect to the $x_1$-direction. 
The boundary $\Gamma:=\partial D$ is supposed to be given by a non-self-intersecting Lipschitz 
curve which is bounded in $x_2$-direction and $2\pi$-periodic with respect to $x_1$. 
Therefore, in this paper we exclude the case of Figure~1 (b) but refer to 
Section~\ref{sec:N}. Let $\tilde{D}$ be a local 
perturbation of $D$ in the way that $\Gamma\setminus\tilde{\Gamma}$ and 
$\tilde{\Gamma}\setminus\Gamma$ are bounded where $\tilde{\Gamma}=\partial \tilde{D}$ is 
the perturbed boundary which is also assumed to be a non-self-intersecting curve. Suppose 
that $\tilde{D}$ is filled by a homogeneous and isotropic medium and that $\tilde{\Gamma}$ 
is a perfectly reflecting curve of Dirichlet kind. Denote by $f\in L^2(\tilde{D})$ a source 
term of compact support which radiates wave fields at the wavenumber $k>0$.
\medlf
We consider the problem of determining the radiated wave $u\in H^1_{loc}(\tilde{D})
:=\bigl\{w|_{\tilde{D}}:w\in H^1_{loc}(\real^2)\bigr\}$ such that
\begin{equation} \label{eqn1}
\Delta u+k^2u\ =\ -f\mbox{ in }\tilde{D}\,,\quad u=0\mbox{ on }\tilde{\Gamma}\,,
\end{equation}
and complemented by the open waveguide radiation condition explained in the next section. 
Without loss of generality (changing the period of the periodic structure if otherwise) we 
can assume that the perturbations $\Gamma\setminus\tilde{\Gamma}$ and 
$\tilde{\Gamma}\setminus\Gamma$ and also the support of $f$ are contained in the disc 
$\{x\in\real^2:(x_1-\pi)^2+x_2^2<\pi^2\}$. We fix $R>\pi$ \tr{and $h_0>\pi$} throughout this paper and 
use the following notations for $h>\pi$ (see Figure~\ref{f12} (a) and 
Figure~\ref{f3}).
\begin{eqnarray*}
Q_h\ & := & \{x\in D:0<x_1<2\pi,\ x_2<h\}\,,\\
\quad Q_\infty\ &:=&\ \{x\in D:0<x_1<2\pi\}\,, \\
\Gamma_h & := & (0,2\pi)\times\{h\}\,,\quad \\
W_h\ &:=&\ \{x\in D:x_2<h\}\,,\quad  \\
 U_h\ &:=&\ \{x\in D: x_2>h\}\, , \\
C_R & := & \{x\in D:(x_1-\pi)^2+x_2^2=R^2\},\quad \\ 
\Sigma_R\ &:=&\ \{x\in D: (x_1-\pi)^2+x_2^2>R^2\}\,, \\
D_R & := & \{x\in D:(x_1-\pi)^2+x_2^2<R^2\}\,,\quad \\
\tilde{D}_R\ &:=&\ \{x\in \tilde{D}:(x_1-\pi)^2+x_2^2<R^2\}.
\end{eqnarray*}
In the unperturbed setting we introduce the following function 
spaces \footnote{The definitions hold also for $D$ instead of $\tilde{D}$}.
\begin{eqnarray*}
H^1_{loc,0}(\tilde{D}) & := & \bigl\{u\in H^1_{loc}(\tilde{D}):u=0\mbox{ on }
\partial\tilde{D}\bigr\}\,, \\
H^1_\ast(\tilde{D}) & := & \bigl\{u\in H^1_{loc}(\tilde{D}):
u|_{W_h\cap\tilde{D}}\in H^1(W_h\cap\tilde{D})\mbox{ for all }h>h_0,\bigr\}\,, \\
H^1_\ast(\Sigma_R) & := & \biggl\{u\in H^1_{loc}(\Sigma_R):
\begin{array}{l} u|_{W_h\cap\Sigma_R}\in H^1(W_h\cap\Sigma_R)\;\mbox{for all}\;h>h_0,\\
u=0\;\mbox{on}\;\partial\Sigma_R\cap\partial D \end{array}\biggr\}\,, \\
H^1_{\alpha,loc}(D) & := & \bigl\{u\in H^1_{loc}(D):u(\cdot,x_2)\mbox{ is 
$\alpha$-quasi-periodic} \bigr\}\,, \\
H^1_{\alpha,loc,0}(D) & := & \bigl\{u\in H^1_{\alpha,loc}(D): u=0\mbox{ on }\partial D
\bigr\}\,.
\end{eqnarray*}

\subsection{The Open Waveguide Radiation Condition And An Energy Formula}
\label{ss-OWC}

As mentioned in the introduction part the diffracted field will have a decomposition into a 
(guided) propagating part and a radiating part. The loss of exponential decay of the radiating 
part is a consequence of the existence of cut-off values while the propagative wave numbers 
determine the behavior of the guided part along the waveguide. We first recall that a function $\phi\in L^2_{loc}(\real)$ is called $\alpha$-quasi-periodic if 
$\phi(x_1+2\pi)=e^{2\pi\alpha i}\phi(x_1)$ for all $x_1\in\real$.
\begin{definition} \label{d-exceptional}
(i) $\alpha\in[-1/2,1/2]$ is called a \emph{cut-off value} if there exists $\ell\in\ganz$ 
such that $|\alpha+\ell|=k$.
\newline
(ii) $\alpha\in[-1/2,1/2]$ is called a \emph{propagative wave number} if there exists a 
non-trivial $\phi\in H^1_{\alpha,loc,0}(D)$ such that 
\begin{eqngroup}\begin{equation} \label{exc:a}
\Delta\phi + k^2\phi\ =\ 0\text{ in }D\,,
\end{equation}
and $\phi$ satisfies the upward Rayleigh expansion
\begin{equation} \label{exc:b}
\phi(x)\ =\ \sum_{\ell\in\ganz}\phi_\ell\,e^{i(\ell+\alpha)x_1}\,
e^{i\sqrt{k^2-(\ell+\alpha)^2}(x_2-h_0)}\quad\mbox{for }x_2>h_0
\end{equation}\end{eqngroup}
for some $\phi_\ell\in\cmplx$ where the convergence is uniform for $x_2\geq h_0+\eps$ for 
every $\eps>0$. The functions $\phi$ are called guided (or 
propagating or Floquet) modes. 
\end{definition}
In all of the paper, we choose the square root function to be 
holomorphic in the cutted plane $\cmplx\setminus(i\real_{\leq 0})$. In particular, 
$\sqrt{t}=i\sqrt{|t|}$ for $t\in\real_{<0}$. 
In Definition~\ref{d-exceptional} we restrict the quasi-periodic parameter $\alpha$ to the 
interval $[-1/2, 1/2]$, because an $\alpha$-quasi-periodic function must be also 
$(\alpha+j)$-quasi-periodic for any $j\in \N$. Throughout this paper we make the following 
assumptions.
\begin{assumption} \label{assump1}
Let $|\ell+\alpha|\not=k$ for every propagative wave number $\alpha\in[-1/2,1/2]$ and every 
$\ell\in\ganz$; that is, no cut-off value is a propagative wave number.
\end{assumption}
\smalf
Under Assumption \ref{assump1} it can be shown (see, e.g. \cite{KL18} for the case of a flat 
curve $\Gamma=\Gamma_0$ and an additional index of refraction) that at most a finite number of 
propagative wave numbers exists in the interval $[-1/2,1/2]$. Furthermore, if $\alpha$ is a 
propagative wave number with mode $\phi$ then $-\alpha$ is a propagative wave number with 
mode $\overline{\phi}$. Therefore, we can number the propagative wave numbers in $[-1/2,1/2]$ 
such that they are given by $\{\hat{\alpha}_j:j\in J\}$ where $J\subset\ganz$ is finite and 
symmetric with respect to $0$ and $\hat{\alpha}_{-j}=-\hat{\alpha}_j$ for $j\in J$. 
Furthermore, it is known that (under Assumption~\ref{assump1}) every mode $\phi$ is evanescent; 
that is, exponentially decaying as $x_2$ tends to infinity in $D$; that is, satisfies 
$|\phi(x)|\leq c\,e^{-\delta|x_2|}$ for $x_2\geq h_0$ and some $c,\delta>0$ which are 
independent of $x$. The corresponding space
\begin{equation} \label{X_j}
X_j\ :=\ \bigl\{ \phi\in H^1_{\hat{\alpha}_j,loc,0}(D):u\mbox{ satisfies (\ref{exc:a}) and 
(\ref{exc:b}) for }\alpha=\hat{\alpha}_j \bigr\}
\end{equation}
of modes is finite dimensional with some dimension $m_j>0$. On $X_j$ we define the sesqui-linear 
form $B:X_j\times X_j\to\cmplx$ by 
\begin{equation} \label{e-sesqui}
B(\phi,\psi)\ :=\ -2i\int\limits_{Q_\infty}\frac{\partial\phi}{\partial x_1}\,
\overline{\psi}\,dx\,,\quad\phi,\psi\in X_j\,.
\end{equation}
Note that $B$ is hermitian. We make the assumption that $B$ is non-degenerated on every 
$X_j$; that is,
\begin{assumption} \label{assump2}
For every $j\in J$ and $\psi\in X_j$, $\psi\not=0$, the linear form  
$B(\cdot,\psi):X_j\to\cmplx$ is non-trivial on $X_j$; that is, there exists $\phi\in X_j$ 
with $B(\phi,\psi)\not=0$.
\end{assumption}
The hermitian sesqui-linear form $B$ defines the cones $\{\psi\in X_j:B(\psi,\psi)\gtrless 0\}$ 
of propagating waves traveling to the right and left, respectively. We construct a basis of 
$X_j$ with elements in these cones by taking any inner product $(\cdot,\cdot)_{X_j}$ and 
consider the following eigenvalue problem in $X_j$ for every fixed $j\in J$. Determine 
$\lambda_{\ell,j}\in\real$ and non-trivial $\hat{\phi}_{\ell,j}\in X_j$ with 
\begin{equation} \label{evp}
B(\hat{\phi}_{\ell,j},\psi)\ =\ -2i\int\limits_{Q_\infty}\frac{\partial\hat{\phi}_{\ell,j}}
{\partial x_1}\,\overline{\psi}\,dx\ =\ 
\lambda_{\ell,j}\,\bigl(\hat{\phi}_{\ell,j},\psi\bigr)_{X_j}\quad\mbox{for all }\psi\in X_j
\end{equation}
and $\ell=1,\ldots,m_j$. We normalize the basis such that
$\bigl(\hat{\phi}_{\ell,j},\hat{\phi}_{\ell^\prime,j}\bigr)_{X_j}=\delta_{\ell,\ell^\prime}$ 
for $\ell,\ell^\prime=1,\ldots,m_j$. Then $\lambda_{\ell,j}=B(\hat{\phi}_{\ell,j},
\hat{\phi}_{\ell,j})$ and the function $\psi\in X_j$ in Assumption~\ref{assump2} must 
take the form $\psi=\sum_{\ell=1}^{m_j} c_\ell \hat{\phi}_{\ell,j}$ with $c_\ell\neq 0$ for some
$\ell=1,2,\cdots, m_j$. Choosing $\phi=\hat{\phi}_{\ell,j}$, one deduces $B(\phi,\psi)=
c_\ell\lambda_{\ell,j}$. Hence, the Assumption~\ref{assump2} is equivalent to 
$\lambda_{\ell,j}\not=0$ for all $\ell=1,\ldots,m_j$ and $j\in J$.
\medlf
\begin{remark}
\begin{itemize}
\item[(i)] The set of propagative wave numbers obviously depends on $k\in\R_+$. Analogously, 
one may define $k_\ell=k_\ell(\alpha)$ for $\alpha\in[-1/2,1/2]$ as the wave number if the 
problem \eqref{exc:a} and \eqref{exc:b} admits a non-trivial solution. Since the solutions are 
in $H^1_{\alpha,loc,0}(D)$ the values $k_\ell(\alpha)$ are just eigenvalues of $-\Delta$
with respect to $\alpha$-quasi-periodic boundary conditions on the vertical boundary of 
$Q_\infty$ and homogeneous Dirichlet boundary condition on $\Gamma$. The functions
$\alpha\rightarrow k_\ell(\alpha)$ are well known as the dispersion relations/curves.
Throughout our paper the wavenumber $k$ is fixed. Under the Assumption \ref{assump1}, the set 
$\{(\hat{\alpha}_j(k_0),k_0)\}_{j\in J}$ constitutes the intersection points of the dispersion
curves with the line $k=k_0$ in the $(\alpha,k)$-plane. Assumption~\ref{assump1} implies the 
absence of flat dispersion curves.
\item[(ii)] The eigenvalue problem \eqref{evp} originates from the limiting absorption principle
(LAP) by applying an abstract functional theorem that goes back to \cite{KL18}. We refer to 
\cite{FJ16,KL-MMAS} for detailed discussions in justifying the radiation conditions for closed 
full and half-waveguide problems. Note that the choice of the inner product in $X_j$
relies on the way how to perturb the original scattering problem by applying the LAP.
For example, if the LAP is applied to the wavenumber $k$ then $(\phi,\psi)_{X_j}=
\int_{Q_\infty}\phi\overline{\psi}dx$.
\end{itemize}
\end{remark}
 In all 
of the paper we make Assumptions~\ref{assump1} and \ref{assump2} without mentioning this 
always. \trr{The one-dimensional Fourier transform is defined as
$$ (\cF\phi)(\omega)\ :=\ \frac{1}{\sqrt{2\pi}}\int\limits_{-\infty}^\infty\phi(s)\,
e^{-is\omega}\,ds\,,\quad\omega\in\real\,.$$
It can be considered as an unitary operator from $L^2(\real)$ onto itself.}
Now we are able to formulate the radiation condition caused by compactly supported source 
terms, which will also serve as the radiation condition of the Green's function to perturbed 
and unperturbed scattering problems (see Theorem \ref{wps} and Remark \ref{rem:usc}).
\begin{definition} \label{d-RC}
Let $\psi_+,\psi_-\in C^\infty(\real)$ be any functions with $\psi_\pm(x_1)=1$ for 
$\pm x_1\geq\sigma_0$ (for some $\sigma_0>\max\{R,2\pi\}+1$) and $\psi_\pm(x_1)=0$ 
for $\pm x_1\leq\sigma_0-1$. 
\smalf
A solution $u\in H^1_{loc}(\Sigma_R)$ of (\ref{eqn1}) satisfies the open 
waveguide radiation condition with respect to an inner product $(\cdot,\cdot)_{X_j}$ in $X_j$ 
if $u$ has in $\Sigma_R  $ a decomposition into $u=u_{rad}+u_{prop}$ which satisfy the 
following conditions.
\begin{itemize}
\item[(a)] The propagating part $u_{prop}$ has the form 
\begin{equation} \label{u2}
u_{prop}(x)\ =\ \sum_{j\in J}\biggl[\psi_+(x_1)\sum_{\ell:\lambda_{\ell,j}>0}
a_{\ell,j}\,\hat{\phi}_{\ell,j}(x)\ +\ \psi_-(x_1)\sum_{\ell:\lambda_{\ell,j}<0}a_{\ell,j}\,
\hat{\phi}_{\ell,j}(x)\biggr]
\end{equation}
for $x\in\Sigma_R  $ and some $a_{\ell,j}\in\cmplx$. Here, for every $j\in J$ the 
scalars $\lambda_{\ell,j}\in\real$ and $\hat{\phi}_{\ell,j}\in\hat{X}_j$ for 
$\ell=1,\ldots,m_j$ are given by the eigenvalues and corresponding eigenfunctions, 
respectively, of the self adjoint eigenvalue problem (\ref{evp}). Note that by the choice of 
$\psi_\pm$ the propagating part vanishes for $|x_1|<\sigma_0-1$ and is therefore well 
defined in $\Sigma_R  $.
\item[(b)] The radiating part $u_{rad}\in H^1_\ast(\Sigma_R)$ satisfies the generalized 
angular spectrum radiation condition
\begin{equation} \label{angulaa-spectrum-rc}
\int\limits_{-\infty}^\infty\left|\frac{\partial(\cF u_{rad})(\omega,x_2)}{\partial x_2}-
i\sqrt{k^2-\omega^2}\,(\cF u_{rad})(\omega,x_2)\right|^2d\omega\ \longrightarrow\ 0\,,\quad 
x_2\to\infty\,.
\end{equation}
\end{itemize}
\end{definition}
\smalf
The radiation condition \eqref{angulaa-spectrum-rc} can be used to prove well-posedness of the 
Helmholtz equation with a source term which is supported in $x_1$-direction and exponentially 
decays in $x_2$ (see \eqref{eq-source1:a}). It has been shown in \cite{KL18} for the case of 
a half plane problem with an inhomogeneous period layer that the radiation condition of 
Definition \ref{d-RC} for the inner product $(\phi,\psi)_{X_j}=2k\int_{Q_\infty}n\,\phi\,
\overline{\psi}\,dx$ is a consequence of the limiting absorption principle by replacing $k$ 
with $k+i\epsilon$, $\epsilon>0$. In this paper we will not justify this radiation condition, 
although we are sure that this can be done in the same way as \cite{KL18, KL-MMAS}. A second 
motivation of our radiation condition is the following result on the direction of 
the energy flow which will play a central role in the proof of uniqueness.
\begin{lemma}\label{lem:energy} 
Let $u=\sum_{j\in J}\sum_{\ell=1}^{m_j} a_{\ell,j} \hat{\phi}_{\ell,j}$ for some $a_{\ell,j}\in
\cmplx$ and write $q+Q_\infty:=\{x\in D:q<x_1<q+2\pi\}$ for $q\in \R$. Then we have 
$$ 2\,\Im\int\limits_{q+Q_\infty}\overline{u}\,\frac{\partial u}{\partial x_1}\,dx\ =\
\sum_{j\in J}\sum_{\ell=1}^{m_j}\lambda_{\ell,j}\,|a_{\ell,j}|^2\,. $$
\end{lemma}
By Lemma \ref{lem:energy}, the propagating part $u_{prop}$ satisfies the energy formula
$$ 2\,\Im\int\limits_{q+Q_\infty}\overline{u_{prop}}\,\frac{\partial u_{prop}}{\partial x_1}\,
dx\ =\ \left\{\begin{array}{cl}\displaystyle\sum_{j\in J}\sum_{\lambda_{\ell,j}>0}
\lambda_{\ell,j}\,|a_{\ell,j}|^2\,, & q>\sigma_0\,, \\
\displaystyle\sum_{j\in J}\sum_{\lambda_{\ell,j}<0}\lambda_{\ell,j}\,|a_{\ell,j}|^2\,, & 
q<-\sigma_0\,, \end{array}\right. $$
where $\sigma_0>2\pi+1$ is the number specified in Definition \ref{d-RC}.
To prove Lemma \ref{lem:energy}, we have to modify the arguments of \cite{K22} for 
inhomogeneous layered media, because solutions of the Dirichlet and Neumann boundary value 
problems are in $H^1_{loc}(\tilde{D})$ but fail to be in $H^2_{loc}(\tilde{D})$ if 
$\tilde{\Gamma}$ is Lipschitz. For $C^2$-smooth boundaries, the quantity in 
Lemma~\ref{lem:energy} also equals to $4\pi\Im\int_{D\cap\{x_1=\tr{q}\}} \overline{u}\,
\frac{\partial u }{\partial x_1}\,ds$; see \cite[Lemma 6.3]{KL18} and \cite[Lemma 2.6]{K22}.
\medlf
\textbf{Proof of Lemma \ref{lem:energy}.} 
We recall the following form of Green's formula valid in any Lipschitz domain $\Omega$: For 
$u\in H^1(\Omega)$ with $\Delta u\in L^2(\Omega)$ we have
$$ \int\limits_\Omega[\nabla u\cdot\nabla\psi+\psi\,\Delta u]\,dx\ =\ 0\quad\mbox{for all }
\psi\in H^1_0(\Omega)\,. $$
Let $j\in\{1,2\}$. First we show for $\alpha_j$-quasi-periodic solutions $u_j\in 
H^1_{loc}(D)$ of $\Delta u_j+k^2u_j=0$ in $D$ with $u_j=0$ on $\Gamma$ and 
$\alpha_j\in(-1/2,1/2]$ with $\alpha_1\not=\alpha_2$ that
\begin{equation} \label{aux1}
\int\limits_{q+Q_\infty}\biggl[\overline{u_2}\,\frac{\partial u_1}{\partial x_1}-
u_1\,\frac{\partial\overline{u_2}}{\partial x_1}\biggr]\,dx\ =\ 0\,.
\end{equation}
Indeed, defining $\psi(x_1):=1-|x_1-q|/(2\pi)$ and applying Green's theorem in 
$\Omega:=\{x\in D:q-2\pi<x_1<q+2\pi\}$ yields (note that $u_j$ decay exponentially as $x_2$ 
tends to infinity)
\begin{eqnarray*}
0 & = &  \int\limits_\Omega[\nabla u_1\cdot\nabla(\psi\overline{u_2})-k^2(\psi\overline{u_2})\,
u_1]\,dx \\ 
& = & \int\limits_\Omega\psi\,[\nabla u_1\cdot\nabla\overline{u_2}-k^2\overline{u_2}\,u_1]\,dx\ 
+\ \int\limits_\Omega\psi^\prime\,\overline{u_2}\,\frac{\partial u_1}{\partial x_1}\,dx\,.
\end{eqnarray*}
Interchanging the roles of $u_1$ and $\overline{u_2}$ and subtraction yields
\begin{eqnarray*}
0 & = & \int\limits_\Omega\psi^\prime\,\biggl[\overline{u_2}\,\frac{\partial u_1}{\partial x_1} 
- u_1\,\frac{\partial\overline{u_2}}{\partial x_1}\biggr]\,dx \\ 
& = & \int\limits_{q+Q_\infty}\psi^\prime\,\biggl[\overline{u_2}\,\frac{\partial u_1}
{\partial x_1} - u_1\,\frac{\partial\overline{u_2}}{\partial x_1}\biggr]\,dx 
 + \int\limits_{q-2\pi+Q_\infty}\psi^\prime\,\biggl[\overline{u_2}\,\frac{\partial u_1}
{\partial x_1} - u_1\,\frac{\partial\overline{u_2}}{\partial x_1}\biggr]\,dx \\
& = & \frac{1}{2\pi}\int\limits_{q+Q_\infty}\biggl[\overline{u_2}\,\frac{\partial u_1}
{\partial x_1} - u_1\,\frac{\partial\overline{u_2}}{\partial x_1}\biggr]\,dx -
\frac{1}{2\pi}\int\limits_{q-2\pi+Q_\infty}\biggl[\overline{u_2}\,\frac{\partial u_1}
{\partial x_1} - u_1\,\frac{\partial\overline{u_2}}{\partial x_1}\biggr]\,dx \\
& = & \bigl(1-e^{2\pi i(\alpha_2-\alpha_1)}\bigr)\,\frac{1}{2\pi}
\int\limits_{q+Q_\infty}\biggl[\overline{u_2}\,\frac{\partial u_1}{\partial x_1} - 
u_1\,\frac{\partial\overline{u_2}}{\partial x_1}\biggr]\,dx
\end{eqnarray*}
where we used the quasi-periodicity of $u_j$. This yields (\ref{aux1}).
\smalf
Now we rewrite $u$ as
$$ u\ =\ \sum_{j\in J}\sum_{\ell=1}^{m_j}a_{\ell,j}\,\hat{\phi}_{\ell,j}\ =\ 
\sum_{j\in J}u_j\quad\mbox{with}\quad u_j\ :=\ \sum_{\ell=1}^{m_j}a_{\ell,j}\,
\hat{\phi}_{\ell,j}. $$
Then $u_j$ is $\hat{\alpha}_j$-quasi-periodic. Using (\ref{aux1}) and the orthonormalization 
of $\hat{\phi}_{\ell,j}$, we arrive at
\begin{eqnarray*}
& & 2i\Im\int\limits_{q+Q_\infty}\overline{u }\,\frac{\partial u }{\partial x_1}\,
dx\ =\ \int\limits_{q+Q_\infty}\biggl[\overline{u }\,\frac{\partial u}{\partial x_1}-
u\,\frac{\partial\overline{u }}{\partial x_1}\biggr]\,dx \\
& = & \sum_{j\in J}\sum_{j^\prime\in J}\int\limits_{q+Q_\infty}\biggl[\overline{u_j}\,
\frac{\partial u_{j^\prime}}{\partial x_1} -u_{j^\prime}\,\frac{\partial\overline{u_j}}
{\partial x_1}\biggr]\,dx\ =\ \sum_{j\in J}\int\limits_{q+Q_\infty}\biggl[
\overline{u_j}\,\frac{\partial u_j}{\partial x_1} -
u_j\,\frac{\partial\overline{u_j}}{\partial x_1}\biggr]\,dx \\
& = & 2i\Im\sum_{j\in J}\int\limits_{q+Q_\infty}\overline{u_j }\,\frac{\partial u_j}
{\partial x_1}\,dx\ =\ i \Re\biggl[-2i\sum_{j\in J}\int\limits_{Q_\infty}\overline{u_j}\,
\frac{\partial u_j }{\partial x_1}\,dx\biggr]\ =\ i\sum_{j\in J}\sum_{\ell=1}^{m_j}
\lambda_{\ell,j}\,|a_{\ell,j}|^2\,,
\end{eqnarray*}
which proves the lemma. \qed
\biglf
Below we review a result on the asymptotic behavior of $u_{rad}$ which will be needed in the 
proof of uniqueness. By \eqref{eqn1} and \eqref{u2}, the radiating part $u_{rad}$ to the 
scattering problem satisfies
\begin{eqngroup}\begin{equation} \label{eq-source1:a}
\Delta u_{rad}\ +\ k^2\,u_{rad}\ =\ -f\ -\ \sum_{j\in J}\sum_{\ell=1}^{m_j}
a_{\ell,j}\varphi_{\ell,j}\quad\mbox{in }\tilde{D}\,,\quad u_{rad}=0\mbox{ on }\tilde{\Gamma}\,,
\end{equation}
where
\begin{equation} \label{eq-source1:b}
\varphi_{\ell,j}(x)\ =\ \left\{\begin{array}{cl} 2\,\psi_+^\prime(x_1)\,
\frac{\partial\hat{\phi}_{\ell,j}(x)}{\partial x_1}+\psi_+^{\prime\prime}(x_1)\,
\hat{\phi}_{\ell,j}(x) & \mbox{if }\lambda_{\ell,j}>0\,, \\ 
2\,\psi_-^\prime(x_1)\,\frac{\partial\hat{\phi}_{\ell,j}(x)}{\partial x_1}+
\psi_-^{\prime\prime}(x_1)\,\hat{\phi}_{\ell,j}(x) & \mbox{if }\lambda_{\ell,j}<0\,. 
\end{array}\right.
\end{equation}\end{eqngroup}
We note that $f$ has compact support in $Q_{h_0}$ and $\varphi_{\ell,j}$ vanish for 
$|x_1|\leq\sigma_0-1$ and $|x_1|\geq\sigma_0$, and are evanescent; that is, there exist 
$\hat{c},\delta>0$ with $|\varphi_{\ell,j}(x)|\leq \hat{c}\exp(-\delta x_2)$ for all 
$x_2\geq h_0$. Furthermore, $u_{rad}$ satisfies the generalized angular spectrum radiation 
condition~(\ref{angulaa-spectrum-rc}). In \cite{K22} the following result has been 
shown.\footnote{These properties are consequences of the differential equation and radiation 
condition above the line $x_2=h_0$ solely and are therefore independent of the differential 
equation or boundary condition below this line.}
\begin{lemma} \label{l-decay1}
Let Assumptions~\ref{assump1} and \ref{assump2} hold, and let $u\in H_{loc}^1(D)$ be a 
solution of (\ref{eqn1}) satisfying the radiation condition of Definition~\ref{d-RC}. 
Then the radiating part $u_{rad}$ satisfies a stronger form of the radiation 
condition~(\ref{angulaa-spectrum-rc}), namely,
\begin{equation} \label{Fouriea-rc-strong}
\left|\frac{\partial(\cF u_{rad})(\omega,x_2)}{\partial x_2}-
i\sqrt{k^2-\omega^2}\,(\cF u_{rad})(\omega,x_2)\right| \leq\ 
\frac{c}{\delta+\sqrt{|\omega^2-k^2|}}\,e^{-\delta x_2}
\end{equation}
for almost all $\omega\in\real$ and $x_2>h_0$ where $c>0$ is independent of $\omega$ and $x$.
\smalf
Furthermore, there exists $c>0$ with
\begin{equation} \label{est-rho}
\bigl|u_{rad}(x)\bigr|\ +\ \bigl|\nabla u_{rad}(x)\bigr|\ \leq\ c\,(1+|x_2|)\,\rho(x_1)
\end{equation}
for all $x\in\tilde{D}$ with $x_2\geq h_0+1$, where $\rho\in L^2(\real)\cap L^\infty(\real)$ 
is given by
\begin{equation} \label{rho}
\rho(x_1)\ :=\ \int\limits_\real\frac{|u_{rad}(y_1,h_0)|}{(1+|x_1-y_1|)^{3/2}}\,dy_1\ +\ 
\frac{1}{1+|x_1|^{3/2}}\,,\quad x_1\in\real\,.
\end{equation}
\end{lemma}

\trr{
\subsection{A Modified Open Waveguide Radiation Condition}
In this subsection we propose another open waveguide radiation condition that is equivalent to 
the Def. \ref{d-RC}. We first define the half-plane Sommerfeld radiation condition used 
in \cite{HWR, K22}. Introduce the weighted Sobolev space $H_\rho^1(\Omega)$  by
\ben
H_\rho^1(\Omega)\ :=\ \bigl\{u:(|1+|x_1|^2)^{\rho/2}
u\in H^1(\Omega)\bigr\}\,,\quad \rho\in\R\,.
\enn
\begin{definition}\label{src} 
A function $v\in C^\infty(U_{h_0}\cap\Sigma_R  )$ satisfies the Sommerfeld radiation 
condition in $U_{h_0}\cap\Sigma_R  $ if $v\in H^1_\rho(W_h\cap\Sigma_R  )$ for all 
$h>h_0$ and all $\rho<0$ and
\be\label{SRC}
\sup\limits_{x\in C_a\cap U_h}|x|^{1/2}\bigl|\frac{\partial v(x)}{\partial r}-ik v(x)\bigr|\ 
\rightarrow\ 0\,,\quad a\rightarrow\infty\,,\qquad 
\sup\limits_{x\in U_h}|x|^{1/2}|v(x)|\ <\ \infty
\en
for all $h>h_0$ where $r=|x|$. 
\end{definition}
\begin{remark}
Since $\Phi(x,y)=\cO(|x|^{-1/2})$ and $\frac{\partial\Phi(x,y)}{\partial r}-ik
\Phi(x,y)=\cO(|x|^{-3/2})$ as $r=|x|\rightarrow\infty$, it holds that $\Phi(\cdot,y)\in 
H^1_\rho(W_h\cap\Sigma_R  )$ for all $\rho<0$ if $R>|y_1-\pi|$. Hence, the above Sommerfeld 
radiation condition covers two-dimensional point source waves, but excludes plane waves and 
surface (evanescent) waves, which do not decay along the horizontal direction.
\end{remark}
If $\Gamma$ is a Lipschitz function,
it was shown in \cite{HWR} that the scattered field caused by a point source source must satisfy the above Sommerfeld radiation condition. However, the total field (i.e., the Green's function to the rough surface scattering problem) satisfies an analogous condition but with the weighted index $\rho<1$ in place of $\rho<0$. Motivated by this fact, we define a modified open waveguide radiation condition by changing the generalized angular spectrum radiation condition of the radiating part of Def. \ref{d-RC}.
\begin{definition} \label{d-MRC}
A solution $u\in H^1_{loc}(\Sigma_R)$ of (\ref{eqn1}) satisfies the modified open 
waveguide radiation condition with respect to an inner product $(\cdot,\cdot)_{X_j}$ in $X_j$ 
if $u$ has  a decomposition into $u=u_{rad}+u_{prop}$ in $\Sigma_R$ where
$u_{prop}$ satisfied the same condition specified as in Def. \ref{d-RC} (a) and 
$u_{rad}$ fulfills the Sommerfeld radiation condition of Def. \ref{src} but with the index $\rho<1$. 
\end{definition}
Below we prove the equivalence of the two open waveguide radiation conditions.
\begin{theorem}
The open waveguide radiation condition of Def. \ref{d-RC} and the modified one given by Def. \ref{d-MRC} are equivalent.
\end{theorem}
\textbf{Proof.}
Write $u=u_{rad}+u_{prop}$ where 
$u_{rad}\in H^1_\ast(\tilde{D})$ denotes the radiating part and $u_{prop}$ the propagating part. First we suppose that $u_{rad}$ fulfills the generalized angular spectrum radiation condition \eqref{angulaa-spectrum-rc}.
By arguing analogously to \cite[Theorem~6.2]{K22} for compact source terms, one can show the 
asymptotics $u_{rad}(x)=O(|x_1|^{-3/2})$ as $|x_1|\rightarrow \infty$ in $W_h$. This gives 
$u_{rad}\in H^1_\rho(W_h\cap \Sigma_R)$ for all $h>h_0$, $\rho<1$ and proves the modified open  waveguide radiation condition of Definition~\ref{d-MRC}; see \cite[Section 6]{K22} for details.
\smalf
Now it remains to justify the generalized angular spectrum radiation condition of $u_{rad}$, under the assumption that $u_{rad}$ satisfies  the Sommerfeld radiation condition of Def. \ref{src} but with the index $\rho<1$. Since $u_{rad}|_{\Gamma_{h_0}}\in H^{1/2}_{\rho}(\R)$ for all $1/2<\rho<1$,
we recall from \cite[Lemma A.2, Appendix]{HWR} (see also \cite{K22}) that the function 
$$
v(x)=2\int_{\Gamma_{h_0}}\frac{\partial G(x, y)}{\partial y_2} u_{rad}(y)\,ds(y),\quad x_2>h_0,
$$
satisfies the homogeneous Helmholtz equation together the Sommerfeld radiation conditions \ref{SRC} and the boundary value $v=u_{rad}$ on $x_2=h_0$. Hence, the function $w:=u_{rad}-v$ satisfies \ref{SRC} in $x_2>h_0$ and the boundary value problem
$$
\Delta w+k^2w=\varphi\quad\mbox{in}\quad x_2>h_0,\qquad w=0\quad\mbox{on}\quad x_2=h_0,
$$
where $\varphi$ is given by the right hand side of \eqref{eq-source1:a}.This implies that $w$ can be represented as 
$$
w(x)=\int_{-\sigma_0}^{\sigma_0}\int_{h_0}^\infty [G(x,y)-G(x, y^*)]\,\varphi(y)\,dy_2dy_1,\quad x_2>h_0,
$$
with $y^*:=(y_1-2h_0-y_2)^\top$. Now, following the proof of \cite[Lemma 7.1]{K22} one can show that $w$ satisfies the stronger form \eqref{Fouriea-rc-strong} of the radiation 
condition~(\ref{angulaa-spectrum-rc}).This proves the generalized angular spectrum radiation condition of $u_{rad}$.
\qed
}
\smalf
We would like to extend the Sommerfeld radiation condition up to the boundary $\Gamma$. 
However, since $\Gamma$ is only Lipschitz, in general the derivatives $\partial v/\partial r$ 
do not exist up to the boundary. We can, however, define a weaker form which models the 
integral form of the Sommerfeld radiation condition as follows. The connection between these 
two radiation conditions will be described in Lemma \ref{coro4.1}. 
\begin{definition}\label{src-i} 
Let $a_j$ be a sequence in $\real$ such that $a_j\to\infty$ and $\tilde{D}_{a_j}$ are Lipschitz 
domains. A solution $v\in H^1_{loc}(\Sigma_R)$ 
satisfies the Sommerfeld radiation condition in integral form if 
$$ \biggl\Vert\frac{\partial v}{\partial r}-ik v\biggr\Vert_{H^{-1/2}(C_{a_j})}\ 
\longrightarrow\ 0\,,\quad j\to\infty\,, $$
where $r=|x|$.
\end{definition}

\begin{lemma}\label{coro4.1} If $v$ satisfies the Sommerfeld radiation condition of 
Definition~\ref{src} with the index $\rho\geq 0$, then v also fulfills the integral form of 
the radiation condition defined by Definition~\ref{src-i}.
\end{lemma}
\textbf{Proof.} Without loss of generality we suppose that $a_{j+1}-a_j\geq 1$ for all 
$j\in\N$. Let $h_0>0$ be the number specified in Definition \ref{src}. We set $h:=h_0+1$ and 
choose $\psi\in C^\infty(\real^2)$ such that $\psi(x)=0$ for $x\in U_h$ and $\psi(x)=1$ for 
$x\notin U_{h-\eps}$. We decompose $C_{a_j}$ into $C_{a_j}=\bigl(C_{a_j}\cap U_h\bigr)\cup
\bigl(C_{a_j}\setminus U_h\bigr)$. Then
\be\nonumber
\biggl\Vert\frac{\partial v}{\partial r}-ik v\biggr\Vert_{H^{-1/2}(C_{a_j})} & = &
\biggl\Vert\psi\biggl(\frac{\partial v}{\partial r}-ik v\biggr)\biggr\Vert_{H^{-1/2}(C_{a_j})}\ 
+\ \biggl\Vert(1-\psi)\biggl(\frac{\partial v}{\partial r}-ik v\biggr)
\biggr\Vert_{H^{-1/2}(C_{a_j})} \\ \nonumber
& \leq & \biggl\Vert\psi\frac{\partial v}{\partial r}\biggr\Vert_{H^{-1/2}(C_{a_j})}\ +\
k\bigl\Vert\psi v\bigr\Vert_{H^{-1/2}(C_{a_j})} \\ \label{e1}
& & +\ \biggl\Vert(1-\psi)\biggl(\frac{\partial v}{\partial r}-ik v\biggr)
\biggr\Vert_{H^{-1/2}(C_{a_j})}.
\en
The last integral converges to zero because, by the Sommerfeld radiation condition of 
(\ref{SRC}),
\ben
\biggl\Vert(1-\psi)\biggl(\frac{\partial v}{\partial r}-ik v\biggr)
\biggr\Vert_{H^{-1/2}(C_{a_j})}\ \leq\ C\,\biggl\Vert\frac{\partial v}{\partial r}-ik v
\biggr\Vert_{L^{2}(C_{a_j}\cap U_{h-\epsilon})}\ \longrightarrow\ 0
\enn
as $j\rightarrow\infty$.  It remains to discuss the first two integrals on the right hand side 
of \eqref{e1}. Let $E_j:H^{1/2}_0(C_{a_j})\to H^1(D_{a_j+1}\setminus\overline{D_{a_j}})$ be
extension operators which are uniformly bounded with respect to $j$. In fact, given 
$\varphi\in H^{1/2}_0(C_{a_j})$ we define $E_j\varphi=w_j$ in $D _{a_{j+1}}\setminus
\overline{D _{a_j}}$ where $w_j$ is the unique solution to the boundary value problem
\ben
& & \Delta w_j=0\quad\mbox{in}\quad D _{a_{j+1}}\setminus\overline{D _{a_j}},\\
& & w_j=\varphi\quad\mbox{on}\quad C_{a_j},\qquad w_j=0\quad\mbox{on}\quad
\partial(D_{a_{j+1}}\setminus\overline{D_{a_j}})\setminus C_{a_j}.
\enn
The norm of such an extension operator depends only on the Lipschitz constants of 
$D_{a_{j+1}}\setminus\overline{D_{a_j}}$, which are uniformly bounded in $j$. Then we have
\begin{eqnarray*}
& & \biggl\Vert\psi\,\frac{\partial v}{\partial r}\biggr\Vert_{H^{-1/2}(C_{a_j})}\ +\
k\bigl\Vert\psi v\bigr\Vert_{H^{-1/2}(C_{a_j})} \\
& \leq & \biggl\Vert\frac{\partial(\psi v)}{\partial r}\biggr\Vert_{H^{-1/2}(C_{a_j})}\ +\
\biggl\Vert v\,\frac{\partial\psi}{\partial r}\biggr\Vert_{H^{-1/2}(C_{a_j})}\ +\
k\bigl\Vert\psi v\bigr\Vert_{H^{-1/2}(C_{a_j})} \\
& \leq & \sup_{\Vert\varphi\Vert_{H^{1/2}_0(C_{a_j})}=1}\langle\partial_r(\psi v),\varphi
\rangle\ +\ c\,\bigl\Vert v\bigr\Vert_{H^{1/2}(C_{a_j}\setminus U_h)} \\
& = & \sup_{\Vert\varphi\Vert_{H^{1/2}_0(C_{a_j})}=1}\int\limits_{D_{a_{j+1}}\setminus
\overline{D _{a_j}}}\bigl[\nabla(\psi v)\cdot\nabla\overline{E_j\varphi}-
k^2 \psi v\,\overline{E_j\varphi}\bigr]\,dx\ +\ 
c\,\bigl\Vert v\bigr\Vert_{H^{1/2}(C_{a_j}\setminus U_h)} \\
& \leq & c\,\Vert v\Vert_{H^1(Z_j)}
\end{eqnarray*}
where $Z_j=\{x\in D_{a_{j+1}}:|x|>a_j\,,\ x_2<h\}$ and $c>0$ is independent of $j$. Simple 
estimates show that $Z_j$ is contained in the set $\{x\in D:a_j-\eps<x_1<a_{j+1},\ x_2<h\}$. 
From $v\in H^1(W_h)$ we conclude that $\Vert v\Vert_{H^1(Z_j)}$ tends to zero.  \qed

\subsection{Uniqueness Of Solutions Of The Perturbed And Unperturbed Problems}
\label{ss-unique}
 
First we show that the propagating part $u_{prop}$ \trr{of the open waveguide radiation condition \ref{d-RC} has to vanish, if $f=0$}.

\begin{theorem} \label{t-unique-u2}
Let $u\in H^1_{loc,0}(\tilde{D})$ be a solution of $\Delta u+k^2u=0$ in $\tilde{D}$ 
satisfying the open waveguide radiation condition of Definition \ref{d-RC}. Then $u_{prop}$ 
vanishes; that is, all the coefficients $a_{\ell,j}$ vanish.
\end{theorem}
\smalf
\textbf{Proof.} Choose $\psi_N\in C^\infty(\real)$ and $\varphi_H\in C^\infty(\real)$ 
with $\psi_N(x_1)=1$ for $|x_1|\leq N$ and $\psi_N(x_1)=0$ for $|x_1|\geq N+1$ and 
$\varphi_H(x_2)=0$ for $x_2\geq H+1$ and $\varphi_H(x_2)=1$ for $x_2\leq H$.
\newline
For $N>\sigma_0+1$ and $H>h_0+1$ we define the regions $D_{N,H}:=\{x\in\tilde{D}:|x_1|<N,\ 
x_2<H\}$ and $W^-_{N,H}:=\{x\in\tilde{D}:-N-1<x_1<-N,\ x_2<H\}$ and $W^+_{N,H}:=
\{x\in\tilde{D}:N<x_1<N+1,\ x_2<H\}$ and the horizontal line segments 
$\Gamma_{N,H}:=(-N,N)\times\{H\}$. We apply Green's theorem in 
$D_{N+1,H+1}$ to $v(x):=\psi_N(x_1)\,u(x)$ and $\overline{v(x)}\,\tr{\varphi_H}(x_2)$. First we note 
that $\Delta v\in L^2(D_{N+1,H+1})$ because $\Delta u=-k^2u$. Furthermore $v\varphi\in
H^1_0(D_{N+1,H+1})$, therefore,
\begin{eqnarray*}
0 & = & \int\limits_{D_{N+1,H+1}}\bigl[\nabla v\cdot\nabla(\overline{v}\,\tr{\varphi_H})+
(\overline{v}\,\tr{\varphi_H})\,\Delta v\bigr]\,dx \\ 
& = & \int\limits_{D_{N+1,H}}\bigl[|\nabla v|^2+\overline{v}\,\Delta v\bigr]\,dx\ 
+ \int\limits_{D_{N+1,H+1}\setminus D_{N+1,H}}\bigl[\nabla v\cdot\nabla(\overline{v}\,\tr{\varphi_H})+
(\overline{v}\,\tr{\varphi_H})\,\Delta v\bigr]\,dx \\
& = & \int\limits_{D_{N+1,H}}\bigl[|\nabla v|^2+\overline{v}\,\Delta v\bigr]\,dx\ 
- \int\limits_{\Gamma_{N+1,H}}\overline{v}\,\frac{\partial v}{\partial x_2}\,ds
\end{eqnarray*}
where we applied the classical Green's theorem in the rectangle $(-N-1,N+1)\times(H,H+1)$ 
to the second integral for the smooth function $v$. Therefore,
\begin{eqnarray*}
&&\int\limits_{\Gamma_{N+1,H}}\psi_N^2\,\overline{u}\,\frac{\partial u}{\partial x_2}\,ds
 =  \int\limits_{\Gamma_{N+1,H}}\overline{v}\,\frac{\partial v}{\partial x_2}\,ds\ =\ 
\int\limits_{D_{N+1,H}}\bigl[\bigl|\nabla v\bigr|^2+\overline{v}\,\Delta v\bigr]\,dx \\ 
& = & \int\limits_{D_{N,H}}\bigl[\bigl|\nabla u\bigr|^2+\overline{u}\,\Delta u\bigr]\,dx 
 +\ \int\limits_{W^+_{N,H}}\bigl[\bigl|\nabla v\bigr|^2+\overline{v}\,\Delta v\bigr]\,dx\ 
+\ \int\limits_{W^-_{N,H}}\bigl[\bigl|\nabla v\bigr|^2+\overline{v}\,\Delta v\bigr]\,dx\,;
\end{eqnarray*}
that is, with $\Delta u=-k^2u$,
\begin{equation} \label{aux3}
\Im\int\limits_{\Gamma_{N+1,H}}\psi_N^2\,\overline{u}\,\frac{\partial u}{\partial x_2}\,ds\
=\ \Im\int\limits_{W^+_{N,H}}\bigl[\bigl|\nabla v\bigr|^2+\overline{v}\,\Delta v\bigr]\,dx\ +\ 
\Im\int\limits_{W^-_{N,H}}\bigl[\bigl|\nabla v\bigr|^2+\overline{v}\,\Delta v\bigr]\,dx \,.
\end{equation}
The decomposition $u=u_{rad}+u_{prop}$ yields four terms in each of the integrals of 
(\ref{aux3}).
\newline 
(a) First, we look at the two integrals on the right hand side of (\ref{aux3}). We 
define $v^{(1)}:=\psi_Nu_{rad}$ and $v^{(2)}=\psi_Nu_{prop}$ and estimate the terms
$$ a^\pm_{N,H}(j,\ell)\ :=\ \int\limits_{W^\pm_{N,H}}\bigl[\nabla\overline{v^{(j)}}\cdot
\nabla v^{(\ell)}+\overline{v^{(j)}}\,\Delta v^{(\ell)}\bigr]\,dx $$
for $j,\ell\in\{1,2\}$. Then $|a^\pm_{N,H}(1,1)|$, $|a^\pm_{N,H}(1,2)|$, and 
$|a^\pm_{N,H}(2,1)|$ are estimated as in the proof of \cite[Theorem~2.2]{K22}:
$$ |a^\pm_{N,H}(1,1)|\ \leq\ c\,\gamma_{N,H}\,,\qquad
|a^\pm_{N,H}(1,2)|\ +\ |a^\pm_{N,H}(2,1)|\ \leq\ c\,\sqrt{\gamma_{N,H}} $$
with
\begin{equation} \label{gamma_N}
\gamma_{N,H}\ :=\ \Vert u_{rad}\Vert_{H^1(Q_N)}^2\ +\ H^3\int\limits_{N<|x_1|<N+1}
\rho(x_1)^2\,dx_1 
\end{equation}
and $Q_N:=\{x\in\tilde{D}:N<|x_1|<N+1,\ x_2<h_0+1\}$.
\newline
For $a^\pm_{N,H}(2,2)$ we need to argue differently as in the proof of \cite[Theorem~3.2]{K22} 
to avoid the integral over the vertical boundaries of $W_{N,H}^\pm$. We recall that
$$ a^+_{N,H}(2,2)\ =\ \int\limits_{W^+_{N,H}}\bigl[|\nabla(\psi_N\,u_{prop})|^2+
(\psi_N\,\overline{u_{prop}})\,\Delta(\psi_Nu_{prop})\bigr]\,dx $$
and note that $\Delta(\psi_N\,u_{prop})=-k^2\psi_N\,u_{prop}+2\psi_N^\prime\,
\frac{\partial u_{prop}}{\partial x_1}+\psi_N^{\prime\prime}\,u_{prop}$. Therefore,
\begin{eqnarray*}
\Im a^+_{N,H}(2,2) & = & 2\Im\int\limits_{W^+_{N,H}}\psi_N\,\psi_N^\prime\,
\overline{u_{prop}}\,\frac{\partial u_{prop}}{\partial x_1}\,dx\ =\ \Im\int\limits_{W^+_{N,H}}
\frac{d}{dx_1}\psi_N^2\,\overline{u_{prop}}\,\frac{\partial u_{prop}}{\partial x_1}\,dx \\
& = & \Im\int\limits_{W^+_{N,H}}\bigl[\nabla u_{prop}\cdot\nabla(\psi_N^2\,\overline{u_{prop}})
-k^2(\psi_N^2\,\overline{u_{prop}})\,u_{prop}\bigr]\,dx \\
& = & \Im\int\limits_{N+Q_\infty}\bigl[\nabla u_{prop}\cdot\nabla(\psi_N^2\,
\overline{u_{prop}})-k^2(\psi_N^2\,\overline{u_{prop}})\,u_{prop}\bigr]\,dx\ -\ \beta^+_{N,H}
\end{eqnarray*}
where again $N+Q_\infty=\{x\in\tilde{D}:N<x_1<N+2\pi\}$ and 
$$ |\beta^+_{N,H}|\ =\ \biggl|\int\limits_{(N+Q_\infty)\setminus W^+_{N,H}}
\bigl[\nabla u_{prop}\cdot\nabla(\psi_N^2\,\overline{u_{prop}})-
k^2(\psi_N^2\,\overline{u_{prop}})\,u_{prop}\bigr]\,dx\biggr|\ \leq\ c\,e^{-2\delta H}\,. $$
Now we set $\varphi(x_1)=1-(x_1-N)/(2\pi)$ and observe that $\psi_N^2-\varphi$ vanishes for 
$x_1=N$ and $x_1=N+2\pi$. Green's theorem implies that
$$ \int\limits_{N+Q_\infty}\bigl[\nabla u_{prop}\cdot\nabla\bigl((\psi_N^2-\varphi)\,
\overline{u_{prop}}\bigr)-k^2\bigl((\psi_N^2-\varphi)\,\overline{u_{prop}}\bigr)\,
u_{prop}\bigr]\,dx\ =\ 0 $$
and thus
\begin{eqnarray*}
\Im a^+_{N,H}(2,2) & = & \Im\int\limits_{N+Q_\infty}\bigl[\nabla u_{prop}\cdot
\nabla(\varphi\,\overline{u_{prop}})-k^2(\varphi\,\overline{u_{prop}})\,u_{prop}\bigr]\,dx\ -\ 
\beta^+_{N,H} \\
& = & -\frac{1}{2\pi}\Im\int\limits_{N+Q_\infty}\overline{u_{prop}}\,
\frac{\partial u_{prop}}{\partial x_1}\,dx\ -\ \beta^+_{N,H} \\
& = & -\frac{1}{4\pi}\sum_{j\in J}\sum_{\lambda_{\ell,j}>0}\lambda_{\ell,j}\,|a_{\ell,j}|^2\
-\ \beta^+_{N,H}
\end{eqnarray*}
where we used the results of Lemma~\ref{l-decay1} above. The same estimates hold for 
$a_{N,H}^-(j,\ell)$; that is, the integrals over $W^-_{N,H}$. Therefore, we have shown that
\begin{eqnarray}
& & \Im\int\limits_{W^+_{N,H}}\bigl[\bigl|\nabla v\bigr|^2+\overline{v}\,\Delta v\bigr]\,dx\ +\ 
\Im\int\limits_{W^-_{N,H}}\bigl[\bigl|\nabla v\bigr|^2+\overline{v}\,\Delta v\bigr]\,dx 
\label{A1} \\ 
& \leq & -\frac{1}{4\pi}\sum_{j\in J}\sum_{\lambda_{\ell,j}>0}\lambda_{\ell,j}\,|a_{\ell,j}|^2\ 
+\ \frac{1}{4\pi}\sum_{j\in J}\sum_{\lambda_{\ell,j}<0}\lambda_{\ell,j}\,|a_{\ell,j}|^2 \ +\
c\,e^{-2\delta H}\ +\ c\,[\gamma_{N,H}+\sqrt{\gamma_{N,H}}]\,. \nonumber
\end{eqnarray}
(b)
Now we look at the left hand side of (\ref{aux3}). The line integrals are outside of 
the layer $W_{h_0}$. Their estimates in \cite{K22} (proof of Theorem~3.2) are 
independent of the equation or boundary condition below the line $x_2=h_0$. In 
\cite{K22} we have shown the existence of sequences $(N_m)$ and $(H_m)$ converging to 
infinity such that $\gamma_{N_m,H_m}\to 0$ and
\begin{equation} \label{est1}
\limsup_{m\to\infty}\biggl[\Im\int\limits_{\Gamma_{N_m+1,H_m}}\psi_{N_m}^2\,\overline{u}\,
\frac{\partial u}{\partial x_2}\,ds\biggr]\ \geq\ 0\,.
\end{equation}
From (\ref{A1}) we conclude that
\begin{eqnarray*}
& & \limsup\limits_{m\to\infty}\biggl[
\Im\int\limits_{W^+_{N_m,H_m}}\bigl[\bigl|\nabla v\bigr|^2+\overline{v}\,\Delta v\bigr]\,dx\ 
+\ \Im\int\limits_{W^-_{N_m,H_m}}\bigl[\bigl|\nabla v\bigr|^2+\overline{v}\,\Delta v\bigr]\,dx 
\biggr] \\
& \leq & -\frac{1}{4\pi}\sum_{j\in J}\sum_{\lambda_{\ell,j}>0}\lambda_{\ell,j}\,|a_{\ell,j}|^2\ 
+\ \frac{1}{4\pi}\sum_{j\in J}\sum_{\lambda_{\ell,j}<0}\lambda_{\ell,j}\,|a_{\ell,j}|^2\,.
\end{eqnarray*}
Combining this estimate with (\ref{est1}) and (\ref{aux3}) yields that $a_{\ell,j}=0$ for all 
$\ell$ and $j$. \qed
\smalf
\trr{
Below we sketch another proof based on the modified open waveguide radiation condition.
\begin{theorem}\label{Th2.15} 
Let $u\in H^1_{loc,0}(\tilde{D})$ be a solution of $\Delta u+k^2u=0$ in $\tilde{D}$ 
satisfying the modified open waveguide radiation condition of Definition \ref{d-MRC}. Then $u_{prop}$ 
vanishes.\end{theorem}
\textbf{Proof.}
Choose $a>0$ and suppose without loss of generality that $\tilde{D}_a$ is a Lipschitz domain. Applying Green's formula for $u$ to $\tilde{D}_a$ gives 
$$
\int_{\tilde{D}_a} |\nabla u|^2-k^2 |u|^2\,dx+\int_{C_a} \partial_\nu u\overline{u}\,ds=0.
$$ Here we have used the Dirichlet boundary condition on $\partial \tilde{D}$. Taking the imaginary part and using $u=u_{rad}+u_{prop}$ yields 
$$
0=\Im \int_{C_a} \partial_\nu u\overline{u}\,ds=\Im \int_{C_a} [\partial_\nu u_{rad}\overline{u}_{rad}+\partial_\nu u_{rad}\overline{u}_{prop}+\partial_\nu u_{prop}\overline{u}_{rad}+ \partial_\nu u_{prop}\overline{u}_{prop} ]\,ds
$$
Recalling the Sommerfeld radiation condition of $u_{rad}$, one can show that the integrals involving the term $u_{rad}$ on the right hand side all vanish as $a\rightarrow\infty$; see the proof of \cite[Theorem 3.1]{HK23} for details. Therefore, one arrives at
$$
0=\lim_{a\rightarrow\infty}\Im \int_{C_a} \partial_\nu u_{prop}\,\overline{u}_{prop}\,ds= \lim_{a\rightarrow\infty}\Im \int_{\tilde{D}_a} (\Delta+k^2) u_{prop}\,\overline{u}_{prop}\,dx,  $$
because $u_{prop}$ vanishes on $\partial \tilde{D}$.
Rcalling the representation of $u_{prop}$ (see \eqref{u2}) and using the definition of $\psi^\pm$, one deduces that (see e.g., \cite[Lemma 2.3]{HK23})
\begin{eqnarray}\label{prop}
\int_{\tilde{D}_a} (\Delta+k^2) u_{prop}\,\overline{u}_{prop}\,dx 
= \left(\int\limits_{\gamma_a^{+}}- \int\limits_{\gamma_a^{-}}\right)\overline{u_{prop}}\,\frac{\partial u_{prop}}{\partial x_1}\,
ds\ +\int\limits_{S_a^+\cup S_a^{-}} \overline{u_{prop}}\,\frac{\partial u_{prop}}{\partial \nu}\,
ds \end{eqnarray}
where $\gamma_a^{\pm}=\tilde{D}_a\cap \{x_1=\pm \sigma_0\}$ and 
$S_a^{\pm}=\{x\in\tilde{D}_a: |x|=a, \sigma_0-1<\pm x_1<\sigma_0\}$.
Note that here have supposed that $\{x\in \tilde{D}_a:\sigma_0-1<\pm x_1<\sigma_0\}$ are both Lipschitz domains by the choice of $\sigma_0>\max\{R, 2\pi\}+1$ and that the integrals over $\gamma_a^\pm$ are understood in the dual form between $H^{-1/2}(\gamma_a^\pm)$ and $H^{1/2}_0(\gamma_a^\pm)$.
The second term on the right hand side of \eqref{prop} tends to zero as $a\rightarrow\infty$, due to the exponential decay of $u_{prop}$ as $x_2\rightarrow\infty$. For the imaginary part of the first term, we have the limit (see \cite[Lemma 2.6]{K22})
\ben
\Im \int\limits_{\tilde{D}\cap \{x_1=\pm \sigma_0\}} \overline{u_{prop}}\,\frac{\partial u_{prop}}{\partial x_1}\,
ds&=&
\frac{1}{2\pi}\Im\int\limits_{\tilde{D}\cap \{\pm \sigma_0+Q_\infty\}} \overline{u_{prop}}\,\frac{\partial u_{prop}}{\partial x_1}dx \\
&=&
 \frac{1}{4\pi} \displaystyle\sum_{j\in J}\sum_{\pm \lambda_{\ell,j}>0}
\lambda_{\ell,j}\,|a_{\ell,j}|^2.
\enn 
In this step we have used the fact $\pm\sigma_0+Q_\infty$ are Lipschitz domains and Lemma \ref{lem:energy}.
Finally, taking the imaginary part in \eqref{prop} and letting  $a\rightarrow\infty$, we obtain $a_{l,j}=0$ for all $j\in J$ and $l=1,2,\cdots,m_j$. This proves $u_{prop}\equiv 0$.
\qed
}
\biglf
Having proved the unique determination of the propagating part, we can show uniqueness of 
solutions of the unperturbed and perturbed boundary value problems following almost the same 
lines in the proof of \cite[Theorem 3.3]{K22}. We omit the proof of Theorem \ref{th:uni} below.
\begin{theorem}\label{th:uni}
Let the Assumptions \ref{assump1} and \ref{assump2} hold. 
\begin{itemize} 
\item[(i)] Let $u\in H^1_{loc,0}(D)$ be a solution of $\Delta u+k^2u=0$ in $D$ satisfying the 
open waveguide radiation condition of Definition \ref{d-RC}. Then $u\equiv 0$. 
\item[(ii)] In the perturbed case we have $u\equiv 0$, if there are no bound states to the 
problem \eqref{eqn1}, that is, any solution $u\in H^1_0(\tilde{D})$ of $\Delta u+k^2u=0$ in 
$\tilde{D}$ must vanish identically. 
\end{itemize}
\end{theorem}
\medlf
In the remaining part we suppose that there are no bound states for the perturbed scattering 
problem, so that uniqueness always holds true by Theorem \ref{th:uni} (ii). Note that this 
assumption can be removed, if the domain $\tilde{D}$ fulfills the following condition (see 
\cite{SM05}):
\begin{equation}\label{GC}
(x_1,x_2)\in\tilde{D}\quad\Longrightarrow\quad(x_1, x_2+s)\in\tilde{D}\quad\mbox{for all }
s>0\,.
\end{equation}
Obviously, the geometrical condition \eqref{GC} can be fulfilled  if the boundary 
$\tilde{\Gamma}$ is given by the graph of some continuous function. But then also the 
existence of \trr{guided} modes is excluded.

\section{Construction of the Dirichlet-to-Neumann (DtN) operator}

For simplicity we suppose that there is an open arc of the form 
$C_R:=\{x\in D:|x_1-\pi|^2+|x_2|^2=R^2\}$ for some $R>\pi$ such that the domain $D_R$ is 
Lipschitz (otherwise we can replace $C_R$ by an open curve with a slightly different shape). 
This implies that the perturbed defect $\Gamma\setminus\tilde{\G}$ always lies below $C_R$.
We refer to Figure~\ref{f3} for a typical situation.

To reduce the scattering problem to a bounded domain, we need Sobolev spaces defined on an open arc. Define the Sobolev spaces (see \cite{M10})
\ben
H^{1/2}_0(C_R) & := & \bigl\{f\in H^{1/2}(\partial D_R):f=0\mbox{ on }
\partial D\setminus C_R\bigr\}\,, \\
H^{1/2}(C_R) & := & \bigl\{f|_{C_R}:f\in H^{1/2}(\partial D_R)\bigr\}\,.
\enn
An important property of $H^{1/2}_0(C_R)$ is that the zero extension of $u$ to 
$\partial D_R$ belongs to $H^{1/2}(\partial D_R)$. We remark that in the previous definitions 
the closed boundary $\partial D_R $ can be replaced by other closed boundaries. If 
$u\in H^1(D_R)$ with $u=0$ on $\Gamma\cap\overline{D_R }$, then we have the traces 
$u|_{C_R}\in H_0^{1/2}(C_R)$. The spaces $\tr{H^{1/2}_0(C_R)}$ and $H^{-1/2}(C_R)$ are 
(anti-linear) dual spaces in the sense that
\ben
\langle\phi,\psi\rangle_{H^{1/2}_0(C_R),H^{-1/2}(C_R)}\ =\
\langle\tilde{\phi},\psi\rangle_{H^{1/2}(\partial D_R),H^{-1/2}(D_R)}\,,
\enn 
where $\tilde{\phi}$ denotes the zero extension of $\phi$ to $\partial D_R $. We further remark
that for any $a>R$ there exists a bounded extension operator $E$ from $H^{1/2}_0(C_R)$ 
into $H^1_0(D_a)$. Indeed, extending $\psi\in H^{1/2}_0(C_R)$ by zero in 
$\partial D_R \cap\Gamma$ we observe that this extension is in $H^{1/2}(\partial {D}_R)$. By 
well known results for Lipschitz domains there exists a bounded extension operator $E_1$ 
from $H^{1/2}(\partial D_R)$ into $H^1(D_R)$. In the same way one extends 
$\psi\in H^{1/2}_0(C_R)$ by zero in $\partial D_a\setminus\overline{D_R }$ and constructs an 
extension operator $E_2$ from $H^{1/2}(\partial(D_a\setminus D_R))$ into 
$H^1(D_a\setminus D_R)$ with zero boundary values for $|x|=a$.
\medlf
Below we recall the definition of the Floquet-Bloch transform to be used later.
\begin{definition} For $g\in C_0^\infty(\R)$, the Floquet-Bloch transform $F$ is defined by
$$ (Fg)(x_1,\alpha)\ :=\ \sum_{n\in\Z} g(x_1+2\pi n)\,e^{-i2\pi n \alpha}\,,\qquad 
x_1\in\R,\ \alpha\in[-1/2, 1/2]\,. $$
The Floquet-Bloch transform $F$ extends to an unitary operator from $L^2(\R)$ to 
$L^2((0, 2\pi)\times(-1/2,1/2))$. If $g$ depends on two variables $x_1$ and $x_2$ then the 
symbol $F$ means the Floquet-Bloch transform with respect to $x_1$.
\end{definition}
In the next subsection we prepare several auxiliary results before constructing the DtN 
operator. 

\subsection{Existence Results For Some Unperturbed Problems}
\label{ss-aux}

The first result is well known and a simple application of the Theorem of Riesz.
Define the weighted Sobolev spaces $H^s_{(\rho)}=\bigl\{u\in H^s_0(D):w_\rho u\in H^s(D)
\bigr\}$ where $w_\rho(x)=e^{\rho|x|}$ for $\rho\geq 0$.
\begin{theorem} \label{t-1}
Let $\varphi\in H^{-1/2}(C_R)$ and $\rho\in(0,1)$. Then there exists a unique solution
$v\in H^1_0(D)$ of
\begin{equation} \label{eqn-t1}
\int\limits_D[\nabla v\cdot\nabla\overline{\psi}+v\,\overline{\psi}]\,dx\ =\ 
\int\limits_{C_R}\varphi\,\overline{\psi}\,ds\quad\mbox{for all }\psi\in H^1_0(D)\,.
\end{equation}
Note that we have written the dual form $\langle\varphi,\psi\rangle$ on the right hand side as 
integral. Here we need that the trace $\psi|_{C_R}\in H^{1/2}_0(C_R)$. Furthermore, 
$v\in H^1_{(\rho)}(D)$ and $\varphi\mapsto v$ is bounded from 
$H^{-1/2}(C_R)$ into $H^1_{(\rho)}(D)$ and even compact from $H^{-1/2}(C_R)$ into 
$L^2_{(\rho^\prime)}(D)$ for all $\rho^\prime<\rho$.
\end{theorem}
\textbf{Proof:} The left hand side is just the inner product in $H^1(D)$, and the 
right hand side is estimated by 
\ben
\bigl|\int_{C_R}\varphi\,\overline{\psi}\,ds\bigr|\ \leq\ \Vert\varphi\Vert_{H^{-1/2}(C_R)}
\,\Vert\psi\Vert_{H^{1/2}_0(C_R)}\ \leq\ c\Vert\varphi\Vert_{H^{-1/2}(C_R)}\,
\Vert\psi\Vert_{H^1(D)}.
\enn 
Therefore, Riesz's theorem implies uniqueness and existence of a solution in $H^1_0(D)$. Set 
$\tilde{v}=w_\rho v$ and $\tilde{\psi}=\frac{1}{w_\rho}\psi$. Then $\nabla\psi=
\nabla w_\rho\tilde{\psi}+w_\rho\nabla\tilde{\psi}$ and $\nabla v=-\frac{\tilde{v}}{w_\rho^2}
\nabla w_\rho+\frac{1}{w_\rho}\nabla\tilde{v}$. Substituting this into the variational equation 
yields
$$ \int\limits_D\biggl[\nabla\tilde{v}\cdot\nabla\overline{\tilde{\psi}}+
\overline{\tilde{\psi}}\,\frac{\nabla w_\rho}{w_\rho}\cdot\nabla\tilde{v}-\tilde{v}\,
\frac{\nabla w_\rho}{w_\rho}\cdot\nabla\overline{\tilde{\psi}}-
\frac{|\nabla w_\rho|^2}{w_\rho^2}\,\tilde{v}\,\overline{\tilde{\psi}}+
\tilde{v}\,\overline{\tilde{\psi}}\biggr]\,dx\ =\ \int\limits_{C_R}\varphi\,
\overline{\tilde{\psi}}\,w_\rho\,ds\,. $$
We observe that the left hand side defines a sesqui-linear form on $H^1(D)$ which is 
coercive for $\rho<1$ because $\frac{|\nabla w_\rho|}{w_\rho}=\rho$. The right hand side 
defines again a bounded linear form on $H^1(D)$. Therefore, Lax-Milgram yields existence
and uniqueness. This proves that $v\in H^1_{(\rho)}(D)$ and that $\varphi\mapsto v$ is 
bounded from $H^{-1/2}(C_R)$ into $H^1_{(\rho)}(D)$. 
\smalf 
Finally we show that $H^1_{(\rho)}(D)$ is compactly imbedded in $L^2_{(\rho^\prime)}(D)$
for all $\rho^\prime<\rho$. Let $(v_j)$ be a sequence in $H^1_{(\rho)}(D)$ which converges
weakly to zero. Set again $\tilde{v}_j=w_\rho v_j$. Then $(\tilde{v}_j)$ converges 
weakly to zero in $H^1(D)$ and is thus bounded. Therefore, there exists $c>0$ with 
$\Vert\tilde{v}_j\Vert_{L^2(D)}\leq c$ for all $j$. We estimate for any $a>0$
\begin{eqnarray*}
\int\limits_{D\setminus D_a}w^2_{\rho^\prime}(x)\,|v_j(x)|^2\,dx & = &
\int\limits_{D\setminus D_a}w^2_\rho(x)\,|v_j(x)|^2\,e^{-2(\rho-\rho^\prime)|x|}\,dx \\ 
& \leq & e^{-2(\rho-\rho^\prime)a}\int\limits_D w^2_\rho(x)\,|v_j(x)|^2\,dx\ \leq\ 
c^2\,e^{-2(\rho-\rho^\prime)a}\,.
\end{eqnarray*}
Given $\varepsilon>0$ we choose $a>0$ with $c^2\,e^{-2(\rho-\rho^\prime)a}<\frac{\eps^2}{2}$ 
and keep $r$ fixed. Since $(\tilde{v}_j)$ tends to zero weakly in $H^1(D)$ it tends to zero 
weakly in $H^1(D_a)$. Therefore, $\Vert\tilde{v}_j\Vert_{L^2(D_a)}$ tends to zero and 
thus also $\int_{D_a}w_{\rho^\prime}^2|v_j|^2dx$ because on $D_a$ the norms 
$\Vert w_\eta v\Vert_{L^2(D_a)}$ are all equivalent. Thus, for sufficiently large $j$ the 
term $\int_{D_a}w_{\rho^\prime}^2|v_j|^2dx$ is less than $\frac{\eps^2}{2}$. \qed
\biglf
The proofs of most
existence results for the Helmholtz equation in periodic structures are 
based on the following result for quasi-periodic problems. For a proof we refer to 
\cite[Theorems~4.2, 4.3, and Remark~4.4]{K22} adopted to the present situation.
\begin{theorem} \label{t-exist-qp}
Let Assumptions~\ref{assump1} and \ref{assump2} hold and let $g_\alpha\in L^2(Q_\infty)$ 
for $\alpha\in[-1/2,1/2]$ depend continuously differentiable on $\alpha$ in $[-1/2,1/2]$.
Let there exist $\hat{c}>0$ and $\delta>0$ with $|g_\alpha(x)|+
|\partial g_\alpha(x)/\partial\alpha|\leq\hat{c}e^{-\delta x_2}$ for almost all 
$x\in Q_\infty$ with $x_2>h_0$ and all $\alpha\in[-1/2,1/2]$. Furthermore, let 
$G\in H^{-1}(Q_{h_0})=H^1_0(Q_{h_0})^\ast$ and assume that for any propagative wave number 
$\hat{\alpha}_j\in[-1/2,1/2]$ the orthogonality condition
\begin{equation} \label{ortho-qp}
\langle G,\overline{\hat{\phi}}\rangle\ +\ \int\limits_{Q_\infty}g_{\hat{\alpha}_j}(x)\,
\overline{\hat{\phi}(x)}\,dx\ =\ 0 
\end{equation}
hold for all modes $\hat{\phi}\in X_j$ corresponding to the propagative wave number 
$\hat{\alpha}_j$. Here, $\langle\cdot,\cdot\rangle$ denotes the dual (bi-linear) form.
\smalf
Then for every $\alpha\in[-1/2,1/2]$ there exists an $\alpha-$quasi-periodic solution 
$v_\alpha\in H_{\alpha,loc,0}^1(D)$ of the equation
\begin{equation} \label{source-qp}
\Delta v_\alpha +k^2v_\alpha\ =\ -g_\alpha\ -\ G\quad\mbox{in }Q_\infty
\end{equation}
satisfying the generalized Rayleigh radiation condition 
\begin{equation} \label{Rayleigh-RC}
\sum_{n\in\mathbb{Z}}\biggl|\frac{\partial v_{\alpha,n}(x_2)}{\partial x_2}-
i\sqrt{k^2-(\alpha+n)^2}\,v_{\alpha,n}(x_2)\biggr|^2\ \rightarrow\ 0\,,
\qquad x_2\rightarrow+\infty\,.
\end{equation}
Here, $v_{\alpha,n}(x_2)=\frac{1}{\sqrt{2\pi}}\int_0^{2\pi}v_\alpha(x_1,x_2)
e^{-i(n+\alpha)x_1}dx_1$ are the Fourier coefficients of $v_\alpha(\cdot,x_2)$, and 
(\ref{source-qp}) is understood in the variational sense
$$ \int\limits_{Q_\infty}\bigl[\nabla v_\alpha\cdot\nabla\overline{\psi}-k^2v_\alpha\,
\overline{\psi}\bigr]\,dx\ =\ \langle G,\overline{\psi}\rangle\ +\ \int\limits_{Q_\infty}
g_\alpha\,\overline{\psi}\,dx $$
for all $\psi\in H^1_{\alpha,loc,0}(D)$ which vanish for $x_2>h$ for some $h>h_0$.
\smalf
Furthermore, $v_\alpha$ can be chosen to depend continuously on $\alpha$, and 
for every $h>h_0$ there exists $c_h>0$ with
$$ \Vert v_\alpha\Vert_{H^1(Q_h)}\ \leq\ c_h\bigl[
\sup\limits_{\beta\in[-1/2,1/2]}\Vert g_\beta\Vert_{L^{(1,2)}(Q_\infty)}+
\sup\limits_{\beta\in[-1/2,1/2]}\Vert\partial g_\beta/\partial\beta\Vert_{L^{(1,2)}(Q_\infty)}\ 
+\ \Vert G\Vert_{H^{-1}(Q_{h_0})}\bigr] $$
for all $\alpha\in[-1/2,1/2]$ where we used the notation 
$\Vert\phi\Vert_{L^{(1,2)}(Q_\infty)}:=\Vert\phi\Vert_{L^1(Q_\infty)}+
\Vert\phi\Vert_{L^2(Q_\infty)}$.
\end{theorem}
\biglf
We will apply this result to the following two problems.

Given $\varphi\in H^{-1/2}(C_R)$, consider the problem of determining $u\in H^1_{loc}(D)$ 
such that
\begin{equation} \label{eqn2}
\Delta u+k^2u=0\mbox{ in }D\setminus C_R\,,\quad u=0\mbox{ on }\Gamma\,,\quad
\frac{\partial u_-}{\partial\nu}-\frac{\partial u_+}{\partial\nu}=\varphi\mbox{ on }C_R\,,
\end{equation}
and that $u$ satisfies the open waveguide radiation condition of Definition \ref{d-RC}. 
Here the normal direction $\nu$ is supposed to direct into the exterior $\Sigma_R  $. 
Well-posedness of the variational formulation corresponding to the problem \eqref{eqn2} is 
stated as follows. 
\begin{theorem} \label{t-2}
Let $\varphi\in H^{-1/2}(C_R)$. Then there exists a unique solution $u\in H^1_{loc,0}(D)$ of 
\begin{equation}\label{vaa-eqn2}
\int\limits_D[\nabla u\cdot\nabla\overline{\psi}-k^2u\,\overline{\psi}]\,dx\ =\ 
\int\limits_{C_R}\varphi\,\overline{\psi}\,ds\quad\mbox{for all }\psi\in H^1_{0,c}(D)
\end{equation}
satisfying the open waveguide radiation condition. Here,
$$ H^1_{0,c}(D)\ :=\ \bigl\{\psi\in H^1_0(D):\mbox{there exists $a>0$ with $\psi(x)=0$ 
for }|x|>a\bigr\}\,. $$
Furthermore, the mapping $\varphi\mapsto u|_{C_R}$ is bounded from $H^{-1/2}(C_R)$ into
$H^{1/2}_0(C_R)$.
\end{theorem}
\textbf{Proof:} Uniqueness follows directly from Theorem \ref{t-unique-u2}, part~(i). 
To prove existence, we suppose without loss of generality that $C_R   $ is chosen to lie 
in $Q_{h_0}$ and define the coefficients $a_{\ell,j}$ explicitly as
\begin{equation} \label{a_ell}
a_{\ell,j}\ :=\ \frac{2\pi i}{|\lambda_{\ell,j}|}\int\limits_{C_R}\varphi(x)\,
\overline{\hat{\phi}_{\ell,j}(x)}\,ds(x)\,,\ \ell=1,\ldots,m_j\,,\ j\in J\,.
\end{equation}
Then the propagating part $u_{prop}$ is defined, and the radiating part $u_{rad}$ has to
satisfy
\begin{eqnarray}
& & \int\limits_D[\nabla u_{rad}\cdot\nabla\overline{\psi}-k^2u_{rad}\,\overline{\psi}]\,dx
\label{vaa-u_rad} \\ 
& = & -\int\limits_D[\nabla u_{prop}\cdot\nabla\overline{\psi}-k^2u_{prop}\,\overline{\psi}]\,
dx\ +\ \int\limits_{C_R}\varphi\,\overline{\psi}\,ds \nonumber \\
& = & \sum_{j\in J}\sum_{\ell=1}^{m_j}a_{\ell,j}\int\limits_D\varphi_{\ell,j}\,
\overline{\psi}\,dx\ +\ \int\limits_{C_R}\varphi\,\overline{\psi}\,ds\quad\mbox{for all }
\psi\in H^1_{0,c}(D) \nonumber
\end{eqnarray}
and the generalized angular spectrum radiation condition~(\ref{angulaa-spectrum-rc}). Here,
$\varphi_{\ell,j}$ are given by (\ref{eq-source1:b}). Defining the distribution 
$f_\varphi\in H^{-1}(D)$ as
$$ \langle f_\varphi,\psi\rangle\ :=\ \int\limits_{C_R}\varphi\,\overline{\psi}\,ds\quad
\mbox{for all}\quad\psi\in H^1_0(D)\,, $$
where the right hand side is understood as the duality between $H^{-1/2}(C_R)$ and 
$H^{1/2}_0(C_R)$, we observe that the variational equation represents the differential 
equation
$$ \Delta u_{rad}+k^2u_{rad}\ =\ -\sum_{j\in J}\sum_{\ell=1}^{m_j}a_{\ell,j}\,
\varphi_{\ell,j}\ -\ f_\varphi\,. $$
One now applies Theorem~\ref{t-exist-qp} to $g_\alpha=\sum_{j\in J}\sum_{\ell=1}^{m_j}
a_{\ell,j}\,F\varphi_{\ell,j}$ and $G=f_\varphi$. The orthogonality condition 
(\ref{ortho-qp}) is satisfied by the choice  of the coefficients (\ref{a_ell}) (see 
\cite[Lemma 5.1]{K22}). Therefore, for all $\alpha\in[-1/2,1/2]$ there exists a solution 
$\hat{v}(\cdot,\alpha)\in H^1_{\alpha,loc,0}(D)$ of the $\alpha-$quasi-periodic problems
\begin{eqnarray}
& & \int\limits_{Q_\infty}\bigl[\nabla\hat{v}(\cdot,\alpha)\cdot\nabla\overline{\hat{\psi}}-
k^2\hat{v}(\cdot,\alpha)\,\overline{\hat{\psi}}\bigr]\,dx \label{vaa-u_rad-qp} \\ 
& = & \sum_{j\in J}\sum_{\ell=1}^{m_j}a_{\ell,j}\int\limits_{Q_\infty}
(F\varphi_{\ell,j})(\cdot,\alpha)\,\overline{\hat{\psi}}\,dx\ +\
\int\limits_{C_R}\varphi\,\overline{\hat{\psi}}\,ds \nonumber
\end{eqnarray}
for all $\hat{\psi}\in H^1_\alpha(D)$ which vanish for $x_2>h$ for some $h>h_0$
satisfying the generalized Rayleigh radiation condition (\ref{Rayleigh-RC}). Furthermore, 
$\hat{v}(\cdot,\alpha)$ depends continuously on $\alpha$ and for every $h>h_0$ there exists 
$c_h>0$ with 
\begin{equation} \label{est2}
\Vert\hat{v}(\cdot,\alpha)\Vert_{H^1(Q_h)}\ \leq\
c_h\,\biggl[\sum_{j\in J}\sum_{\ell=1}^{m_j} |a_{\ell,j}|\ +
\sup\limits_{\Vert\psi\Vert_{H^1(Q_\infty)=1}}\biggl|\int\limits_{C_R}\varphi\,
\overline{\psi}\,ds\biggr|\biggr]\ \leq\ c_h\,\Vert\varphi\Vert_{H^{-1/2}(C_R)}\,.
\end{equation}
By the properties of the Floquet-Bloch transform the inverse transform
$$ u_{rad}(x)\ :=\ (F^{-1}\hat{v})(x)\ =\ \int\limits_{-1/2}^{1/2}\hat{v}(x,\alpha)\,d\alpha $$ 
is in $H^1_\ast(D)$. Furthermore, taking $\psi\in C^\infty_0(D)$ we substitute 
$\hat{\psi}:=(F\psi)(\cdot,\alpha)$ into the variational equation~(\ref{vaa-u_rad-qp}) 
and integrate with respect to $\alpha$; that is,
\begin{eqnarray*}
& & \int\limits_{-1/2}^{1/2}\int\limits_{Q_\infty}\bigl[\nabla\hat{v}(x,\alpha)\cdot
\nabla\overline{(F\psi)(x,\alpha)}-k^2\hat{v}(x,\alpha)\,
\overline{(F\psi)(x,\alpha)}\bigr]\,dx\,d\alpha \\ 
& = & \sum_{j\in J}\sum_{\ell=1}^{m_j}a_{\ell,j}\int\limits_{-1/2}^{1/2}\int\limits_{Q_\infty}
(F\varphi_{\ell,j})(x,\alpha)\,\overline{(F\psi)(x,\alpha)}\,dx\,d\alpha +
\int\limits_{-1/2}^{1/2}\int\limits_{C_R}\varphi(x)\,\overline{(F\psi)(x,\alpha)}\,ds(x)\,
d\alpha
\end{eqnarray*}
Noting that $\int_{-1/2}^{1/2}(F\psi)(\cdot,\alpha)\,ds=\psi$ and using the unitarity of the 
Floquet-Bloch transform we observe that this is exactly the equation (\ref{vaa-u_rad}).
Boundedness of $\varphi\rightarrow u|_{C_R}$ is now easily seen from (\ref{est2}) and the 
unitarity of $F$ and the fact that $u_{prop}$ depends explicitly on $\varphi$ through 
$a_{\ell,j}$. \qed
\biglf
A second application is the following result where the source fails to be compactly 
supported.
\begin{theorem} \label{t-3}
Let $f\in L^2_{(\rho)}( {D})$ for some $\rho\in(0,1)$. Then there exists a unique solution 
$w\in H^1_{loc,0}(D)$ of $\Delta w+k^2w=-f$ in $D$ satisfying the open waveguide 
radiation condition. Furthermore, for every $a>R$ the mappings $f\mapsto w|_{D_a}$ and 
$f\mapsto w|_{C_R}$ are bounded from $L^2_{(\rho)}(D)$ into $H^1(D_a)$ and 
$H^{1/2}_0(C_R)$, respectively.
\end{theorem}
\textbf{Proof:} Since $f$ decays exponentially, its Floquet-Bloch transform $Ff$ is well 
defined and continuously differentiable with respect to $\alpha$. Instead of 
(\ref{vaa-u_rad-qp}) we now solve
\begin{eqnarray*}
& & \int\limits_{Q_\infty}\bigl[\nabla\hat{w}(\cdot,\alpha)\cdot\nabla\overline{\hat{\psi}}-
k^2\hat{w}(\cdot,\alpha\,\overline{\hat{\psi}}\bigr]\,dx  \\ 
& = & \sum_{j\in J}\sum_{\ell=1}^{m_j}a_{\ell,j}\int\limits_{Q_\infty}
(F\varphi_{\ell,j})(\cdot,\alpha)\,\overline{\hat{\psi}}\,dx\ +\
\int\limits_{Q_\infty}(Ff)(\cdot,\alpha)\,\overline{\hat{\psi}}\,dx \nonumber
\end{eqnarray*}
for all $\hat{\psi}\in H^1_\alpha(D)$ which vanish for $x_2>h$ for some $h>h_0$
and the generalized Rayleigh radiation condition (\ref{Rayleigh-RC}).

The coefficients $a_{\ell,j}$ have to be chosen as
$$ a_{\ell,j}\ =\ \frac{2\pi i}{|\lambda_{\ell,j}|}\int\limits_{Q_\infty}
(Ff)(x,\hat{\alpha}_j)\,\overline{\hat{\phi}_{\ell,j}(x)}\,dx\,,\quad 
\ell=1,\cdots, m_j\,,\ j\in J\,, $$
so that the right hand side is always orthogonal to the nullspace of the homogeneous 
equation. Using the estimate 
\begin{eqnarray*}
\Vert(Ff)(\cdot,\beta)\Vert^2_{L^2(Q_\infty)} & = & \int\limits_{Q_\infty}
\bigl|(Ff)(x,\beta\bigr|^2dx\ \leq\ \int\limits_{Q_\infty}
\biggl[\sum_{\ell\in\ganz}|f(x_1+2\pi\ell,x_2)|\biggr]^2dx \\
& = & \int\limits_{Q_\infty}\biggl[\sum_{\ell\in\ganz}(1+\ell^2)^{-1/2}\bigl|
(1+\ell^2)^{1/2}f(x_1+2\pi\ell,x_2)\bigr|\biggr]^2dx \\
& \leq & \int\limits_{Q_\infty}\sum_{\ell\in\ganz}\frac{1}{1+\ell^2}\sum_{\ell\in\ganz}
(1+\ell^2)|f(x_1+2\pi\ell,x_2)|^2\,dx \\
& \leq & c\int\limits_D(1+x_1^2)\,|f(x)|^2\,dx\ \leq\ c\,\Vert f\Vert^2_{L^2_{(\rho)}(D)}
\end{eqnarray*}
and analogously for $\Vert(Ff)(\cdot,\beta)\Vert^2_{L^1(Q_\infty)}$ and the derivatives 
with respect to $\beta$ we can repeat the proof of Theorem \ref{t-2}. \qed

\subsection{The DtN Operator}
\label{ss-dtn}

Now we turn to the construction of the Dirichlet-to-Neumann operator on the artificial boundary 
$C_R$. In the remaining part of this paper we make the following assumption.
\begin{assumption}\label{a3}
Assume that $k^2$ is not the Dirichlet eigenvalue of $-\Delta$ in the Lipschitz domain $D_R$ 
and there are no bound states of the Helmholtz equation over the domain $\Sigma_R$; that is, 
if $u\in H^1_0(\Sigma_R)$ solves $\Delta u+k^2u=0$ in $\Sigma_R$, then $u$ must vanish 
identically.
\end{assumption}
As usual, the DtN operator $\Lambda  $ should be defined as follows.
\begin{definition}\label{dtn0} 
The Dirichlet-to-Neumann operator $\Lambda:H_0^{1/2}(C_R)\rightarrow H^{-1/2}(C_R)$ is 
defined by $\Lambda\,g=\partial_\nu u|_{C_R}$ where $u\in H_{loc}^1(\Sigma_R)$ is the 
unique solution to
\begin{equation} \label{dtn}
\Delta u+k^2u=0\mbox{ in }\Sigma_R  \,\quad u=g\mbox{ on }C_R   \,,\quad
u=0\mbox{ on }\Gamma\cap\partial\Sigma_R  \,,
\end{equation} 
which fulfills the open waveguide radiation condition of Definition \ref{d-RC}. Here the unit
normal vector $\nu$ at $C_R   $ is supposed to direct into $\Sigma_R  $. 
\end{definition}
The above definition assumes already the solvability of a boundary value problem in the 
perturbed region $\Sigma_R  $ -- which to show is the purpose of the forthcoming 
Section~\ref{s-pert}. However, the perturbed region $\Sigma_R  $ is a subset of $D$ 
(in contrast to the more general perturbation $\tilde{D}$) which allows the application of an 
integral equation approach with the Dirichlet-Green's function of $D$. Before we explain the 
construction we note that an explicit representation in form of a series can be obtained if 
$\Gamma$ is a straight line parallel to the $x_1$-axis. In this exceptional case the 
propagating part (guided waves) vanishes identically and the radiating part fulfills the 
classical Sommerfeld radiation condition. Consequently, the function $g\in H_0^{1/2}(C_R)$ 
can be expanded into $g(\theta)=\sum_{n\in\N_0}g_n\sin n\theta$ with $\theta\in(0, \pi)$ 
and the DtN operator takes the explicit form
\begin{equation*}
(\Lambda  g)(\theta)\ =\ \sum_{n\in\N_0}g_n\frac{k\,H_n^{(1)\prime}(kR)}{H_n^{(1)}(kR)}\,
\sin n\theta\quad\mbox{for }\theta\in(0,\pi)\,.
\end{equation*}
Here, $H_0^{(1)}$ denotes the Hankel function of the first kind of order zero. In the general 
case that $\Gamma$ is a periodic curve, we will express the field in $\Sigma_R  $ as a 
single layer potential with density $\varphi$ and the Green's function as kernel. As 
usual, $\varphi$ is determined from $g$ by solving an integral equation for the single layer 
boundary operator. We divide our arguments into two steps.
\medlf
\emph{(A) Construction of the single layer boundary operator.} As a motivation we recall 
that for smooth data the single layer boundary operator with the Green's function as kernel is 
given by $S\varphi=u|_{C_R}$ where $u\in H^1_{loc}(D)$ satisfies the transmission 
problem~\eqref{eqn2} and the open waveguide radiation condition. In this way we avoid the 
explicit use of the Green's function. For given $\varphi\in H^{-1/2}(C_R)$, the variational 
form of \eqref{eqn2} is given by \eqref{vaa-eqn2} and has been studied in Theorem~\ref{t-2}.

We take the solution of this transmission problem as the definition of the single layer 
operator, namely $S\varphi:=u|_{C_R}$ where $u\in H^1_{loc,0}(D)$ is the unique solution of
(\ref{vaa-eqn2}) satisfying the open waveguide radiation condition. Then $S$ is bounded from 
$H^{-1/2}(C_R)$ into $H^{1/2}_0(C_R)$ by Theorem~\ref{t-2}. To show the injectivity of $S$, 
we suppose that $S\varphi=0$. Then $u=0$ in $\Sigma_R$ and $u=0$ in $D_R$ by the 
Assumption~\ref{a3} and the uniqueness result of Theorem~\ref{th:uni}. From the 
variational equation (\ref{vaa-eqn2}) we conclude that 
$\int_{C_R}\varphi\,\overline{\psi}\,ds=0$ for all $\psi$ which implies that $\varphi$ 
vanishes. This proves the injectivity of $S$. Next we show that $S$ is boundedly invertible.
\smalf
Let $S_i$ be the operator corresponding to wave number $k=i$. Then, setting $S_i\varphi=
v|_{C_R}$, we get $\langle\varphi,S_i\varphi\rangle=\int_{C_R}\varphi\,\overline{v}\,ds$ by 
the definition of the dual form $\langle\cdot,\cdot\rangle$ (see Theorem \ref{t-1}) 
and $v\in H^1_0(D)$ solves (\ref{eqn-t1}). Setting $\psi=v\,\phi_a$ in (\ref{eqn-t1}) where 
$\phi_a\in C^\infty_0(D)$ satisfies $\phi_a=1$ for $|x|\leq a$ and letting $a$ tend to 
infinity shows that
$$ \langle\varphi,S_i\varphi\rangle\ =\ \int\limits_{C_R}\varphi\,\overline{v}\,ds\ =\ 
\int\limits_D[|\nabla v|^2+|v|^2]\,dx\ =\ \Vert v\Vert^2_{H^1(D)}\,. $$
Next we note that $\Vert\varphi\Vert_{H^{-1/2}(C_R)}=\sup\bigl\{|\langle\varphi,\psi\rangle|:
\Vert\psi\Vert_{H^{1/2}_0(C_R)}\leq 1\bigr\}$. For $\psi\in H^{1/2}_0(C_R)$ with 
$\Vert\psi\Vert_{H^{1/2}_0(C_R)}\leq 1$ we set $\tilde{\psi}=E\psi$ with the extension operator
$E$ from $H^{1/2}_0(C_R)$ into $H^1_0(D_a)$ for some $a>R$ and estimate
$$ |\langle\varphi,\psi\rangle|\ = \biggl|\int\limits_{C_R}\varphi\,\overline{\psi}\,ds\biggr|
= \biggl|\int\limits_ {D}[\nabla v\cdot\nabla\overline{\tilde{\psi}}+
v\,\overline{\tilde{\psi}}]\,dx\biggr| \leq\ c\,\Vert\tilde{\psi}\Vert_{H^1(D)}\,
\Vert v\Vert_{H^1(D)}\ \leq\ c\,\Vert E\Vert\,\Vert v\Vert_{H^1(D)} $$
for $\Vert\psi\Vert_{H^{1/2}_0(C_R)}\leq 1$ and thus $\Vert\varphi\Vert_{H^{-1/2}(C_R)}\leq 
c\,\Vert E\Vert\,\Vert v\Vert_{H^1(D)}$. Combining this with the previous estimate yields
coercivity of $S_i$; that is,
$$ \langle\varphi,S_i\varphi\rangle\ \geq\ \frac{1}{c^2\Vert E\Vert^2}\,
\Vert\varphi\Vert_{H^{-1/2}(C_R)}^2\,. $$
Now we show that $S-S_i$ is compact. We observe that $(S-S_i)\varphi=w|_{C_R}$ where 
$w=u-v\in H^1_{loc}(D)$ satisfies
$$ \Delta w+k^2 w\ =\ -(k^2+1)\,v\mbox{ in }D\,,\quad w=0\mbox{ on }\Gamma\,, $$
and the open waveguide radiation conditions.  Here, $v$ corresponds to the solution of 
\eqref{eqn-t1} with $k=i$ as before. By Theorem~\ref{t-1} we know that $\varphi\mapsto v$ 
is compact from $H^{-1/2}(C_R)$ into $L^2_{(\rho^\prime)}(D)$ for all $\rho^\prime<\rho$. 
Furthermore, by Theorem~\ref{t-3} (for $\rho^\prime$ replacing $\rho$) the mapping 
$(1+k^2)v\mapsto w|_{C_R}$ is bounded from $L^2_{(\rho^\prime)}(D)$ into $H^{1/2}_0(C_R)$. 
Combining this yields compactness of $\varphi\mapsto w|_{C_R}$; that is, compactness of 
$S-S_i$ from $H^{-1/2}(C_R)$ into $H^{1/2}_0(C_R)$.
\medlf
Therefore, the operator equation $S\varphi=g$ can be written as $S_i\varphi+(S-S_i)\varphi=g$. 
This shows that $S$ is a Fredholm operator with index zero. By the Fredholm alternative, the 
injectivity implies the invertibility of $S$.
\medlf
\emph{(B) Construction of the Dirichlet-to-Neumann operator.} Given $g\in H^{1/2}_0(C_R)$ we 
define $\varphi:=S^{-1}g\in H^{-1/2}(C_R)$. Then, by definition, $g=S\varphi=u|_{C_R}$ where 
$u$ satisfies (\ref{vaa-eqn2}); in particular $\Delta u+k^2 u=0$ in $\Sigma_R  $ and $u=0$ 
on $\Gamma\cap \partial \Sigma_R  $, complemented by the open waveguide radiation condition. 
Consequently, the Neumann boundary data can be defined by Green's first formula; that is, the 
DtN operator $\Lambda$ from $H^{1/2}_0(C_R)$ into $H^{-1/2}(C_R)=
\bigl(H^{1/2}_0(C_R)\bigr)^\ast$ can be defined as follows.
\begin{definition} \label{d-DtN}
Let $a>R$ be fixed. Then $\Lambda  :H^{1/2}_0(C_R)\to H^{-1/2}(C_R)=
\bigl(H^{1/2}_0(C_R)\bigr)^\ast$ is defined as
\begin{equation} \label{Lambda} 
\langle\Lambda g,\psi\rangle\ =\ -\int\limits_{D_a\setminus D_R }
\bigl[\nabla u\cdot\nabla\overline{(E\psi)}-k^2u\,\overline{(E\psi)}\bigr]\,dx\,,
\quad\psi\in H^{1/2}_0(C_R)\,,
\end{equation}
where $E:H^{1/2}_0(C_R)\to H^1_0(D_a)$ is again a fixed extension operator and 
$u\in H^1_{loc}(D)$ is the single layer potential with density $\varphi:=S^{-1}g\in 
H^{-1/2}(C_R)$; that is, the unique solution of (\ref{vaa-eqn2}) studied in Theorem~\ref{t-2}.
\end{definition}
We note that the definition is independent of $a>R$ or the choice of the extension operator 
$E$. This follows from Green's identity $\int_{D_a\setminus D_R }
\bigl[\nabla u\cdot\nabla(\overline{\psi_1}-\overline{\psi_2})-k^2u\,
(\overline{\psi_1}-\overline{\psi_2})\bigr]\,dx=0$ for all $\psi_j\in H^1_0(D_a)$ with 
$\psi_1=\psi_2$ on $C_R   $.
\medlf
We finish this section by proving some mapping properties of $\Lambda$.
\begin{lemma}\label{lem:dtn}
The DtN operator $\Lambda:H^{1/2}_0(C_R)\rightarrow H^{-1/2}(C_R)$ is bounded. Moreover, 
the operator $-\Lambda$ can be decomposed into the sum of a coercive operator and a compact 
operator.
\end{lemma}
\textbf{Proof.} By the trace lemma and \eqref{Lambda}, the boundedness of $\Lambda$ follows 
from the estimate 
$$ \Vert\Lambda g\Vert_{H^{-1/2}(C_R)}\ =\ \sup\limits_{\Vert\psi\Vert_{H^{1/2}_0(C_R)}=1}
|\langle\Lambda g,\psi\rangle|\ \leq\ c\,\Vert u\Vert_{H^1(D_a\setminus D_R)}\ \leq\
c\,\Vert g\Vert_{H^{1/2}(C_R)}\,, $$
where we have used the boundedness of the extension operator $E$ and the continuous dependence
of $u$ from $g$. Define $\Lambda_i:H^{1/2}_0(C_R)\rightarrow H^{-1/2}(C_R)$ as the 
DtN operator for the wave number $k=i$; that is,
$$ \langle\Lambda_ig,\psi\rangle\ =\ -\int\limits_{D_a\setminus D_R}
[\nabla v\cdot\nabla\overline{\tilde{\psi}}+v\,\overline{\tilde{\psi}}]\,dx\,,\quad
\psi\in H^{1/2}_0(C_R)\,, $$
where $v\in H^1_0(D)$ solves (\ref{eqn-t1}) for $\varphi:=S_i^{-1}g\in 
H^{-1/2}(C_R)$, and $\tilde{\psi}\in H^1_0(D_a)$ is an extension of $\psi$. The operator 
$-\Lambda_i$ is coercive over $H_0^{1/2}(C_R)$. Indeed, choose $\phi_a\in C^\infty(\real^2)$ 
with $\phi_a=1$ for $|x|<R$ and $\phi_a=0$ for $|x|>a-1$ and set $\tilde{\psi}=v\,\phi_a$ 
for $a>R+1$. Then $\tilde{\psi}\in H^1_0(D_a)$ and thus
$$ -\langle\Lambda_ig,g\rangle\ =\ \Vert v\Vert_{H^1(D_{a-1}\setminus D_R)}^2\ +\ 
\int\limits_{D_a\setminus D_{a-1}}[\nabla v\cdot\nabla\overline{(v\phi_a)}+
\phi_a|v|^2]\,dx\,. $$
Now we let $a$ tend to infinity and use that $v\in H^1(D)$. Therefore,
$$ -\langle\Lambda_ig,g\rangle\ =\ \Vert v\Vert_{H^1(D\setminus D_R)}^2\ \geq\ 
c\,\Vert v\Vert_{H^{1/2}(C_R)}^2\ =\ c\,\Vert g\Vert_{H^{1/2}(C_R)}^2 $$
where we used the boundedness of the trace operator in the inequality. Furthermore, the 
operator $\Lambda-\Lambda_i$ is compact. Indeed, this follows from
$$ \langle(\Lambda-\Lambda_i)\,g,\psi\rangle\ =\ -\int\limits_{D_a\setminus D_R}
\bigl[\nabla(u-v)\cdot\nabla\overline{(E\psi)}-(k^2u+v)\,\overline{E\psi}\bigr]\,dx\,,\quad
\psi\in H^{1/2}_0(C_R)\,, $$
and the compactness of the mapping $g\mapsto(u-v)|_{D_a\setminus D_R}$ from 
$H^{1/2}_0(C_R)$ into $H^1(D_a\setminus D_R)$ (by the same arguments as in the proof of the 
compactness of $S-S_i$) and the boundedness of $g\mapsto k^2u+v$ from 
$H^{1/2}_0(C_R)$ into $H^1(D_a\setminus D_R)$ and the compact embedding of 
$H^1(D_a\setminus D_R)$ into $L^2(D_a\setminus D_R)$. \qed

\section{Existence of Solutions of the Perturbed Problem}
\label{s-pert}

\begin{figure}[htb]
  \centering
  \includegraphics[width=0.65\textwidth]{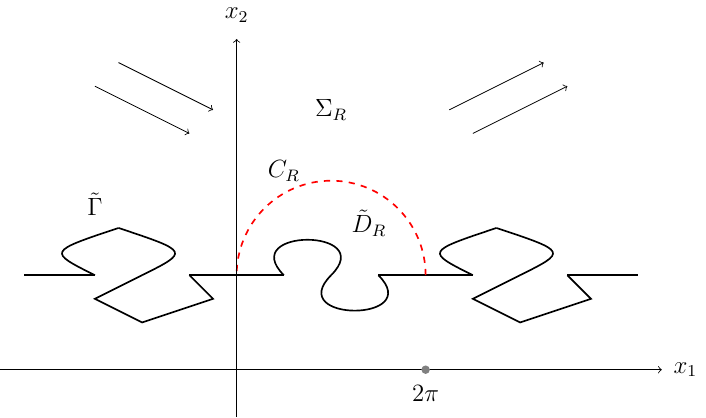}
  \caption{Illustration of wave scattering from perfectly reflecting periodic curves with a 
  local perturbation.}\label{f3}
\end{figure}

In this section we investigate well-posedness of time-harmonic scattering of an incoming wave 
$u^{in}$ from a locally perturbed periodic curve $\tilde{\Gamma}=\partial \tilde{D}$ of 
Dirichlet kind; see Figure \ref{f3}. We consider three kinds of incoming waves: 
\begin{itemize}
\item[(i)] Point source wave: $u^{in}(x):=\Phi(x,y)=\frac{i}{4}H_0^{(1)}(k|x-y|)$ with the 
source position $y\in\tilde{D}$. Without loss of generality we suppose that 
$y\in\tilde{D}_R $. 
\item[(ii)] Plane wave: $u^{in}(x)=e^{ikx\cdot\hat{\theta}}$ where 
$\hat{\theta}=(\sin\theta,-\cos\theta)$ is the incident direction with some incident angle 
$\theta\in(-\pi/2, \pi/2)$. In this case the incoming wave is incident onto $\tilde{\Gamma}$ 
from above, and the parameter $\alpha:=k\sin\theta$ is supposed to be not a propagative wavenumber 
(see Definition \ref{d-exceptional} (ii)). 
\item[(iii)] $u^{in}(x)=\hat{\phi}_{\ell,j}(x)$ is a right (resp. left) going surface wave 
at the propagative wavenumber $\hat{\alpha}_j$ for some $j\in J$ which corresponds to the 
spectral problem~\eqref{evp} with the eigenvalue $\lambda_{\ell,j}>0$ (resp. 
$\lambda_{\ell,j}<0$).
\end{itemize}

We denote by $u^{sc}_{unpert}$ the unperturbed scattered field, defined in $D$, which is 
caused by the unperturbed curve $\Gamma$. In $\Sigma_R  $ the total field $u$ can be decomposed 
into $u=u^{in}+u^{sc}_{unpert}+u^{sc}_{pert}$, and $u^{sc}_{pert}$ can be considered as the 
scattered part induced by the defect. The field $u^{sc}_{pert}$ is supposed to fulfill the open 
waveguide radiation condition of Definition \ref{d-RC} for all of the cases (i), (ii), (iii).

\medlf
Define the spaces 
\begin{eqnarray*}
Y_R & := & \bigl\{v\in H^1(\tilde{D}_R):v=0\mbox{ on }\tilde{\Gamma}\cap\partial\tilde{D}_R   
\bigr\}\,,
\end{eqnarray*}
where $y\in \tilde{D}$. Well-posedness of our 
scattering problems will be stated separately for different incoming waves.

\begin{theorem}[Well-posedness for point source waves]\label{wps}
Let $u^{in}:=\Phi(\cdot,y)$ be an incoming point source wave with $y\in\tilde{D}_R$. Then 
the locally perturbed scattering problem admits a unique solution \cred{$u$}
such that $u-u^{in}\in H^1_{loc}(\tilde{D})$ and $u$ satisfies the \trr{open waveguide radiation 
conditions of Definitions \ref{d-RC} and \ref{d-MRC}}.
\end{theorem}
In this theorem, the total field $u$ is required to satisfy the open waveguide radiation 
condition of  Definition \ref{d-RC}, because $u$ is nothing else but the Green's function of
the perturbed problem. We remark that in $\Sigma_R$ the scattered field
$u-u^{in}$ does not fulfill this radiation condition, since 
$u^{in}=\Phi(\cdot,y)$ does not belong to $H^1(W_h\cap\Sigma_R)$ for any $h>h_0$.

If both $\Gamma$ and $\tilde{\G}$ can be represented as graphs of Lipschitz 
functions, the propagating part $u_{prop}$ vanishes identically (see \cite{SM05,CE10}) and 
thus $u=u_{rad}$. In such a case, it was verified in \cite[Theorem 2.2]{HWR} that 
$u=u(\cdot, y)\in H^1_\rho(W_h\cap\Sigma_R)$ with $R>|y_1-\pi|$ for all $\rho<1$ and that 
$u(\cdot,y)-\Phi(\cdot,y)\in H^1_\rho(W_h)$ for all $\rho<0$. In addition,  both $u$ and 
$u-\Phi$ satisfy the Sommerfeld radiation condition of Definition \ref{src}. The above 
results of Theorem \ref{wps} have generalized those of \cite{HWR} to non-graph curves
where guided (propagating) waves may occur. On the other hand, the technical assumption 
made in \cite[Section 2.3]{HWR} that $\Gamma$ should contain at least one line segment 
in each period  was removed in this paper by constructing a new form of the DtN operator; 
see subsection \ref{ss-dtn}. 
\medlf
\textbf{Proof of Theorem \ref{wps}.} Since the incident field is singular at $y$, 
we transform our scattering problem to an equivalent source problem of the form 
\eqref{eqn1}. Introduce a smooth cut-off function 
$\chi: \mathbb{R}^2\rightarrow \mathbb{R}$ such that $\chi(x)=1$ for $|x-y|<\epsilon/2$ and 
$\chi(x)=0$ for $|x-y|\geq \epsilon$. Here $\epsilon>0$ is chosen to be less than the 
distance between $y$ and $\partial\tilde{D}_R $. We make the ansatz on the total 
field $u$ as
\ben
u(x)\ =\ \chi(x)\,\Phi(x,y)\ +\ v(x,y)\,,\quad x\in\tilde{D},\quad x\neq y\,.
\enn
Then the scattering problem is equivalent of finding $v(\cdot,y)\in 
H_{0, loc}^1(\tilde{D})$ such that
\ben
\left\{\begin{array}{l}\displaystyle \Delta_xv(\cdot,y)+k^2v(\cdot,y)\ =\ -g_y\mbox{ in }
\tilde{D},\quad v(\cdot, y)=0\mbox{ on }\tilde{\Gamma}, \\
\displaystyle \mbox{$v(\cdot, y)$ satisfies the open waveguide radiation condition of 
Definition~\ref{d-RC}} \end{array}\right.
\enn 
with
\ben
g_y\ :=\ \Delta\chi\,\Phi(\cdot,y)+2\nabla\chi\cdot\nabla_x\Phi(\cdot,y)\in L^2(\tilde{D})\,.
\enn
Note that the source term $g_y$ is compactly supported in $\tilde{D}_R$. By the DtN operator 
$\Lambda$, this problem can be reduced to an equivalent boundary value problem over the 
truncated domain $\tilde{D}_R$. Consequently, we get the following variational formulation.
Determine $v\in Y_R$ such that 
\be\label{va:source} 
\int\limits_{\tilde{D}_R}\nabla v\cdot\overline{\nabla \psi}-
k^2v\,\overline{\psi}\,dx\ -\ \int\limits_{C_R}\Lambda  v\,\overline{\psi}ds\ =\
\int\limits_{\tilde{D}_R}g_y\,\overline{\psi}\,ds\quad\mbox{for all}\quad \psi\in Y_R\,.
\en
Here the integral over $C_R$ is understood as the duality between $H^{-1/2}(C_R)$ and 
$H^{1/2}_0(C_R)$. In view of Lemma \ref{lem:dtn}, the sesqui-linear form defined by the left 
hand side of \eqref{va:source} is strongly elliptic, leading to a Fredholm operator with 
index zero over $Y_R$. By Theorem \ref{th:uni} we have uniqueness and thus also 
existence of $v(\cdot,y)\in Y_R$ by the Fredholm alternative. This solution 
can be extended to the exterior $\Sigma_R$ by solving the problem of Theorem \ref{t-1} with 
$\varphi:=S^{-1} (v|_{C_R})\in H^{-1/2}(C_R)$. 
\smalf
Finally we note that $u-u^{in}=v+(\chi-1)\Phi(\cdot,y)\in H^1_{loc}(\tilde{D})$ because $\chi$
vanishes in a neighborhood of $y$, and $u=v$ in $\Sigma_R$ because $\chi$ vanishes in
$\Sigma_R$. This ends the proof. \qed
\biglf
\begin{remark}\label{rem:usc}
If we decompose the field $u$ into $u=u^{in}+u^{sc}_{unpert}+u^{sc}_{pert}$ in $\Sigma_R$ then
we observe that also $u^{sc}_{pert}$ satisfies the open radiation condition because $u$ and 
$u^{in}+u^{sc}_{unpert}$ do, the latter because it is the total field corresponding to the 
unperturbed problem.
\end{remark}
\biglf
We proceed with the scattering problem for plane waves. 
\begin{theorem}[Well-posedness for plane waves]\label{wpw}
Let $\alpha:=k\sin\theta$ be not a propagative wavenumber (see Definition~\ref{d-exceptional} 
(ii)). Then the perturbed scattering problem for a plane wave incidence 
$u^{in}(x)=e^{ikx\cdot\hat{\theta}}$ admits a unique solution 
$u=u^{in}+u^{sc}\in H_{loc,0}^1(\tilde{D})$ such that the scattered part $u^{sc}$ has a 
decomposition in the form $u^{sc}=u^{sc}_{unpert}+u^{sc}_{pert}$ in the region $\Sigma_R  $
where $u^{sc}_{unpert}\in H^1_{\alpha,loc}(D)$ is the scattered field corresponding to the 
unperturbed problem that satisfies the upward Rayleigh expansion~\eqref{exc:b} 
with the quasi-periodic parameter $\alpha=k\sin\theta$. The part $u^{sc}_{pert}\in 
H_{loc}^1(\Sigma_R)$ fulfills \trr{the open waveguide radiation conditions defined by Def. \ref{d-RC} and Def. \ref{d-MRC}}.
\end{theorem}
\textbf{Proof.} In the unperturbed case, uniqueness and existence of the field 
$u^{sc}_{unpert}\in H_{\alpha,loc}^1(D)$ can be justified using standard variational arguments 
in the truncated periodic cell $Q_h$ (for some $h>h_0$) by enforcing the $\alpha-$quasi-periodic
DtN mapping on the artificial boundary $\G_h$. Uniqueness follows from the assumption that
$\alpha=k\sin\theta$ is not a propagative wavenumber, and existence is a consequence of the
Fredholm alternative. 
\smalf
Set $\tilde{u}^{in}=u^{in}+u^{sc}_{unpert}\in H_{\alpha,0}^1(D)$. This field is well 
defined in $\Sigma_R$. We make the ansatz for the perturbed problem in the form 
$u=\tilde{u}^{in}+u^{sc}_{pert}$ in $\Sigma_R$ and $u=u^{in}+u^{sc}$ in $\tilde{D}_R $.
Since $u^{sc}_{pert}\in H^1_{loc}(\Sigma_R)$ is required to fulfill the open waveguide 
radiation condition and $u^{sc}_{pert}=0$ on $\G\setminus\tilde{D}_R$, it has to satisfy
\ben
\partial_\nu u^{sc}_{pert}\bigr|_+\ =\ \Lambda  (u^{sc}_{pert}|_+)\,,\quad\mbox{thus}\quad
\partial_\nu u\bigr|_-\ =\ \partial_\nu\tilde{u}^{in}\bigr|_+\ +\ 
\Lambda  (u|_--\tilde{u}^{in}|_+)\quad \mbox{on }C_R   
\enn
where $|_+$ and $|_-$ denote the traces from $\Sigma_R  $ and $\tilde{D}_R   $, respectively.
Therefore, we have to determine the total field $u\in Y_R$  such that
\be\label{va:plane}
\int\limits_{\tilde{D}_R} \nabla u\cdot\overline{\nabla \psi}-k^2 u\overline{\psi}\,dx\ -\
\int\limits_{C_R} \Lambda u\,\overline{\psi}ds\ =\ \int\limits_{C_R} 
[\partial_\nu\tilde{u}^{in}-\Lambda\tilde{u}^{in}]\,\overline{\psi}\,ds 
\en
for all $\psi\in Y_R$. Application of Lemma \ref{lem:dtn}, Theorem \ref{th:uni} and the 
Fredholm alternative yields the uniqueness and existence of $u\in Y_R$. This also gives the 
scattered field $u^{sc}=u-u^{in}\in H^1(\tilde{D}_R)$ and the trace of the perturbed scattered 
field $g:=u^{sc}_{pert}|_{C_R}=(u^{sc}-u^{sc}_{unpert})|_{-}$ on $C_R$. Finally, 
$u^{sc}_{pert}$ can be extended to $\Sigma_R$ by solving the problem of Theorem \ref{t-1} with 
$\varphi=S^{-1}(g)$. 
\qed
\begin{remark} \label{rem4.8}
Suppose in Theorem \ref{wpw} that $k\sin\theta=\hat{\alpha}_j$ is a propagative 
wavenumber for some fixed $j\in J$. Then it is well known that there exists still a 
$\hat{\alpha}_j-$quasi-periodic solution $u_{unpert,0}=u^{in}+u^{sc}_{unpert}$ of the 
unperturbed problem. However, the solution is not unique, and the general solution is given by
\begin{equation}
u_{unpert}\ =\ u_{unpert,0}\ +\ \sum_{\ell=1}^{m_j}c_\ell\,\hat{\phi}_{\ell, j}\quad
\mbox{in }D
\end{equation} 
where $\hat{\phi}_{\ell,j}\in X_j$ (see \eqref{X_j} and \eqref{evp}) and $c_\ell\in\mathbb{C}$ 
are arbitrary. \tcrr{In our paper \cite{HK23} we derive a new radiation condition based on the limiting absorption principe to prove uniqueness of the unperturbed scattering problem, even if $k\sin\theta$ is a propagative wavenumber}. 
\end{remark}
\biglf
Now we consider the case that $u^{in}=\hat{\phi}_{\ell,j}$ for some $\ell\in\{1,\ldots,m_j\}$ 
and $j\in J$ is an incoming surface wave corresponding to the propagative wavenumber 
$\hat{\alpha}_j$; that is,
\ben\left\{\begin{array}{lll}
\Delta u^{in}+k^2 u^{in}=0\quad\mbox{in}\quad D,\quad
u^{in}=0\quad\mbox{on}\quad \G,\\
\mbox{$u^{in}$ is $\hat{\alpha}_j$-quasi-periodic in $x_1$ and exponentially decays in the 
$x_2$-direction.}
\end{array}\right.
\enn
Since $u^{in}$ vanishes already on $\Gamma$ and satisfies the radiation condition we 
conclude that the variational formulation for $u\in Y_R$ takes the same form as in 
\eqref{va:plane} with $\tilde{u}^{in}=u^{in}$. Analogously to the proof of Theorem \ref{wpw}, 
we obtain
\begin{theorem}[Well-posedness for incoming surface waves]\label{wsw}
Given an incoming surface wave $u^{in}=\hat{\phi}_{\ell,j}$ for some $\ell\in\{1,\ldots,m_j\}$ 
and $j\in J$, the perturbed scattering problem admits a unique solution 
$u=u^{in}+u^{sc}\in H_{loc,0}^1(\tilde{D})$ such that $u^{sc}\in H_{loc}^1(\tilde{D})$ 
fulfills the open waveguide radiation conditions of \trr{Def. \ref{d-RC} and Def. \ref{d-MRC}}. 
\end{theorem}
By Theorem \ref{wsw}, each surface wave $\hat{\phi}_{\ell,j}$ produces a non-trivial scattered 
field to the locally defected problem. Combining Theorems \ref{wpw}, \ref{wsw} and 
Remark~\ref{rem4.8}, we can get a general solution for plane wave incidence when $k\sin\theta$ 
is a progagative wavenumber.

\begin{corollary} 
Let $u^{in}$ be a plane wave and suppose that $k\sin\theta=\hat{\alpha}_j$ is a propagative 
wavenumber for some fixed $j\in J$. The general solution of the perturbed scattering problem for
plane wave incidence takes the form
\begin{equation}
u\ =\ u_{unpert,0}\ +\ u^{sc}_{pert}\ +\ \sum_{\ell=1}^{m_j}c_\ell\,\hat{\phi}_{\ell,j}\
+\ \sum_{\ell=1}^{m_j}c_\ell\, u^{sc}_\ell\quad\mbox{in }\Sigma_R\,.
\end{equation} 
Here, $u_{pert}^{sc}$ is the open waveguide radiation solution determined in Theorem \ref{wpw} 
excited by the incoming reference wave $\tilde{u}^{in}= u_{unpert,0}=u^{in}+u^{sc}_{unpert}$, 
and $u^{sc}_\ell$ is the scattered field specified in Theorem \ref{wsw} with 
$u^{in}:=\hat{\phi}_{\ell,j}$.
\end{corollary}

\section{Scattering by Neumann curves and by periodically arrayed obstacles}\label{sec:N}

With slight changes our solvability results presented in Section \ref{s-pert} carry over to 
periodic and locally perturbed periodic curves of Neumann kind. Below we only remark the 
necessary modifications. 
\medlf
In the Neumann case, $\alpha\in[-1/2,1/2]$ is called a \emph{propagative wave number} if there 
exists a non-trivial $\phi\in H^1_{\alpha,loc}(D)$ such that 
\ben
\Delta\phi + k^2\phi\ =\ 0\text{ in }D\,,\quad\frac{\partial\phi}{\partial\nu}=0\mbox{ on }
\Gamma\,,
\enn
and $\phi$ satisfies the Rayleigh expansion \eqref{exc:b}. 
Here $\nu$ denotes the normal direction at $\Gamma$ pointing into $D$.
Under the Assumption \ref{assump1}, one can still prove that there exist at most a finite 
number of propagative wavenumvers in the interval $[-1/2, 1/2]$. The finite dimensional 
eigenspace $X_j$ can be defined similarly to \eqref{X_j} but with the Neumann boundary 
condition on $\Gamma$. The definition of the space $H^1_\ast(\Sigma_R)$ should be replaced by
\begin{equation*} 
H^1_\ast(\Sigma_R)\ :=\ \bigl\{u\in H^1_{loc}(\Sigma_R):\partial_\nu u=0\mbox{ on }
\Gamma\cap\partial\Sigma_R,\ u\in H^1(W_h\cap\Sigma_R)\mbox{ for all }h>h_0\bigr\}\,.
\end{equation*}

In this case, a bound state of the perturbed scattering problem is defined as a solution 
$u\in H^1(\tilde{D})$ to the Helmholtz equation $(\Delta +k^2)u=0$ in $\tilde{D}$ satisfying 
the Neumann boundary condition $\partial_\nu u=0$ on $\tilde{\Gamma}$. Assuming that there 
are no bound states in $\tilde{D}$, one can prove uniqueness to the perturbed scattering 
problem analogously to Theorem \ref{th:uni}. To construct the DtN operator, we consider the  
problem of determining $u\in H^1_{loc}(D)$ such that, for $\phi\in H_0^{-1/2}(C_R)$,
\be\label{eqn3}
\Delta u+k^2u=0\mbox{ in }D\setminus C_R\,,\quad\partial_\nu u=0\mbox{ on }\Gamma\,,
\quad\frac{\partial u_-}{\partial\nu}-\frac{\partial u_+}{\partial\nu}=\varphi\mbox{ on }C_R\,,
\en
and that $u$ satisfies the open waveguide radiation condition in $\Sigma_R$. The variational 
form of this transmission problem is to determine $u\in H^1_{loc}(D)$ such that
\be\label{eqn4} 
\int\limits_ {D}[\nabla u\cdot\nabla\overline{\psi}-k^2u\,\overline{\psi}]\,dx\ =\ 
\int\limits_{C_R}\varphi\,\overline{\psi}\,ds\quad\mbox{for all }\psi\in H^1_c(D)\,, 
\en
together with the open waveguide radiation condition. Here, 
\ben
H^1_c(D)\ :=\ \{\psi\in H^1(D):\mbox{there exists $a>0$ with $\phi(x)=0$ for all 
$|x|>a$}\}.
\enn
Note that the right hand side of \eqref{eqn4} is understood as the duality between 
$H_0^{-1/2}(C_R)$ and $H^{1/2}(C_R)$. The mapping $S\varphi=u|_{C_R}$ defines the single layer 
operator under the Neumann boundary condition. Choose the open arc $C_R$ such that the 
mixed boundary value problem
\ben
(\Delta+k^2)u=0\mbox{ in }D_R\,,\quad u=0\mbox{ on }C_R\,,\quad 
\partial_\nu u=0\mbox{ on }\partial D_R\setminus C_R\,,
\enn
admits the trivial solution only. We make the assumption that every solution 
$u\in H^1(\Sigma_R)$ to the exterior boundary value problem
\ben
(\Delta+k^2)u=0\mbox{ in }\Sigma_R\,,\quad u=0\mbox{ on }C_R\,,\quad
\partial_\nu u=0\mbox{ on }\Gamma\cap\partial\Sigma_R\,,
\enn 
must vanish identically, that is, there are no bound states to this special perturbation 
problem. The previous two conditions ensure that the single layer operator 
$S:H^{-1/2}_0(C_R)\rightarrow H^{1/2}(C_R)$ is injective and boundedly invertible. The DtN 
operator $\Lambda$ from $H^{1/2}(C_R)$ into $H^{-1/2}_0(C_R)$ takes the explicit form
\ben 
\langle\Lambda g,\psi\rangle\ =\ -\int\limits_{D_a\setminus D_R}
[\nabla u\cdot\nabla\overline{\tilde{\psi}}-k^2u\,\overline{\tilde{\psi}}]\,dx\,,\quad
\psi\in H^{1/2}(C_R)\,, 
\enn
where $\tilde{\psi}=E\psi$ is a bounded extension operator from $H^{1/2}(C_R)$ to 
$H^1_0(D_a\setminus D_R)$ for some $a>R$. Here $u$ is the single layer potential with density 
$\varphi:=S^{-1}g\in H_0^{-1/2}(C_R)$; that is, the open waveguide radiation solution to the 
boundary value problem \eqref{eqn3}. Mapping properties of $\Lambda$ can be proved in the 
same way as Lemma \ref{lem:dtn}. Finally, well-posedness results for scattering of point 
source waves, plane waves and surface waves from locally perturbed Neumann curves can be 
verified in the same manner as in the proofs of Theorems \ref{wps}, \ref{wpw} and \ref{wsw},
\biglf
\begin{figure}[htb]
\centering
\includegraphics[width=0.65\textwidth]{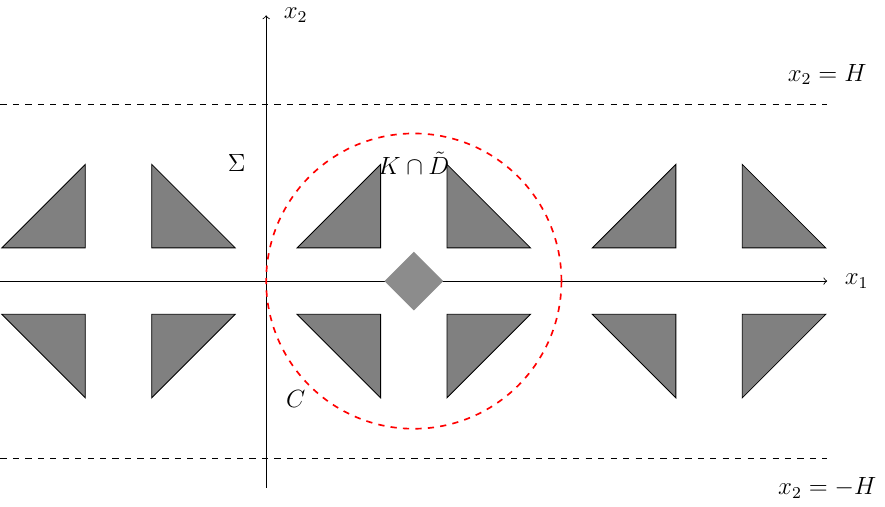}
\caption{Illustration of the artificial boundary $C:=\partial K\subset D$ (in this case 
a circle) on which the DtN operator $\Lambda$ (see Definition \ref{dtn0}) is defined for 
scattering by periodically arrayed obstacles with a local defect. }\label{f4}
\end{figure}

Let us now consider the TE and TM polarizations of time-harmonic electromagnetic scattering 
by periodically arrayed obstacles. Define the boundary conditions $\mathcal{B}u:=u$ in the TE 
case and $\mathcal{B}u:=\partial_\nu u$ in the TM case. 
Let $\Omega\subset\real\times(-H,H)$ be a  domain which is $2\pi-$periodic with respect to 
$x_1$ such that the exterior $D:=\real^2\setminus\overline{\Omega}$ is connected. 
Then $\alpha\in[-1/2,1/2]$ is called a \emph{propagative wave number} if there exists a 
non-trivial $\phi\in H^1_{\alpha,loc}(D)$ such that 
\begin{eqngroup}\ben
\Delta\phi + k^2\phi\ =\ 0\text{ in }D\,,\quad \mathcal{B}\phi=0\mbox{ on }\partial D\,,
\enn
and $\phi$ satisfies the Rayleigh expansions
\ben
\phi(x)\ =\ \sum_{\ell\in\ganz}\phi_\ell^\pm\,e^{i(\ell+\alpha)x_1}\,
e^{\pm i\sqrt{k^2-(\ell+\alpha)^2}(x_2\mp H)}\quad\mbox{for }x_2\gtrless\pm H
\enn\end{eqngroup}
for some $\phi^\pm_\ell\in\cmplx$. Then the spaces $X_j$ of modes and their basis 
$\{\hat{\phi}_{\ell,j}:\ell=1,\ldots,m_j\}$ are defined as in (\ref{X_j})--(\ref{evp}). 
Furthermore, let $\Omega$ be locally defected such that the periodic domain $D$ is replaced by 
a perturbed connected domain $\tilde{D}$. We assume that there exists a bounded Lipschitz 
domain $K$ which contains the defect $(D\setminus\tilde{D})\cup(\tilde{D}\setminus D)$ and  
such that $C:=\partial K$ is contained in $D$. Defining $\Sigma:=D\setminus\overline{K}$ and 
the Sobolev space
\begin{equation*} 
H^1_\ast(\Sigma)\ :=\ \bigl\{u\in H^1_{loc}(\Sigma):\mathcal B u=0\mbox{ on }
\partial D\cap\partial\Sigma\,,\ u\in H^1(W_h\cap\Sigma)\mbox{ for all }h>H\bigr\}\,,
\end{equation*}
where now $W_h:=\real\times(-h,h)$. Then the radiation conditions of Definitions~\ref{d-RC},
\ref{src}, and \ref{src-i} carry over. $C$, $K$, and $\Sigma$ correspond to $C_R$, $D_R$, and 
$\Sigma_R$, respectively. A situation where $C$ can be chosen as a circle $C_R$ is sketched 
in Figure~\ref{f4}. The Dirichlet-to-Neumann operator $\Lambda$ is again defined as 
$\Lambda g=\partial_\nu u|_C$ where $u\in H^1_{loc}(\Sigma)$ is the unique solution of 
\be\label{dtn-obstacle}
\Delta u+k^2u=0\mbox{ in }\Sigma\,,\quad u=g\mbox{ on }C\,,\quad
\mathcal{B}u=0\mbox{ on }\partial D\cap\partial\Sigma\,,
\en
together with the open waveguide radiation condition. We remark that the domain and range space 
of $\Lambda$ relies on the boundary condition under consideration. With proper assumptions on 
the domain $K$, one can construct an invertible single layer operator $S\varphi=u|_C$, 
where $u\in H^1_{loc}(D)$ is the radiating solution of the transmission problem
$$ \Delta u+k^2u=0\mbox{ in }D\setminus C\,,\quad\mathcal{B}u=0\mbox{ on }\partial D\,,\quad
\frac{\partial u_-}{\partial\nu}-\frac{\partial u_+}{\partial\nu}=\varphi\mbox{ on }C\,. $$
Then one can define the DtN operator via Green's formula, analogously to the scattering by 
Dirichlet and Neumann curves. The well-posedness results of Section \ref{s-pert} can be 
justified in the same manner. 

\begin{remark}
Exact boundary conditions (DtN maps) were also constructed for wave propagating in a closed 
periodic waveguide \cite{JLF2006} and in a photonic crystal \cite{FJ09} containing a local 
perturbation. In \tr{comparision} with \cite{FJ09}, the DtN map defined by \eqref{dtn-obstacle} 
applies to \tr{artificial} boundary curves of arbitrary shape (although circular curves are used 
in this paper) and the medium is periodic in one direction only. The exact boundary condition 
of \cite{FJ09} is defined along square-shaped artificial boundaires, and the medium  is 
periodic in two directions. In this paper the DtN map relies heavily on the open waveguide 
radiation condition of Def. \ref{d-RC}.
\end{remark}

\section*{Acknowledgements}
The first author (G.H.) acknowledges the hospitality of the Institute for Applied and Numerical 
Mathematics, Karlsruhe Institute of Technology and the support of Alexander von 
Humboldt-Stiftung.  The second author (A.K.) gratefully acknowledges the financial 
support by the Deutsche Forschungsgemeinschaft (DFG, German Research Foundation) -- 
Project-ID 258734477 -- SFB 1173.

\end{document}